\documentclass[reqno,11pt]{amsart}
\usepackage[top=2.0cm,bottom=2.0cm,left=3cm,right=3cm]{geometry}
\usepackage{amsthm,amsmath,amssymb,dsfont}
\usepackage{mathrsfs,amsfonts,functan,extarrows,mathtools}

\usepackage{xcolor}
\usepackage{cite}
\usepackage{indentfirst, latexsym, amssymb, enumerate,amsmath,graphicx}
\usepackage{float}
\usepackage{relsize}
\usepackage{marginnote}
\usepackage{stmaryrd}
\usepackage{esint}
\usepackage{graphicx}
\usepackage{bm}

\usepackage{cite}
\allowdisplaybreaks[4]
\theoremstyle{plain}
\newtheorem{thm}{Theorem}[section]

\newtheorem{lem}[thm]{Lemma}
\newtheorem{prop}[thm]{Proposition}
\newtheorem{rem}{Remark}[section]

\numberwithin{equation}{section}

\begin{document}

\title[The incompressible  Vlasov-MHD model]
{Global well-posedness and decay rates of strong solutions to the incompressible Vlasov-MHD system
 }

\author[f.-C. Li]{fucai Li}
\address{School of Mathematics, Nanjing University, Nanjing
 210093, P. R. China}
\email{fli@nju.edu.cn}

\author[J.-K. Ni]{Jinkai Ni$^*$} \thanks{$^*$\! Corresponding author}
\address{School of Mathematics, Nanjing University, Nanjing
 210093, P. R. China}
\email{602023210006@smail.nju.edu.cn}

\author[M. Wu]{Man Wu}
\address{Department of Mathematics, Nanjing Audit University, Nanjing
 211815, P. R. China}
\email{manwu@nau.edu.cn}

\begin{abstract}
In this paper, we study the global well-posedness and decay rates of 
strong solutions to  an incompressible  Vlasov-MHD  model arising in  magnetized plasmas.  
This model is consist of the Vlasov equation and the incompressible magnetohydrodynamic equations which interacts   together via the  Lorentz forces.  It is readily to verify that it has two equilibria $(\bar f,\bar u,\bar B)=(0,0,0)$ and $( \tilde f,\tilde u,\tilde B)=(M,0,0)$,
 where $M$ is the global maxwellian.  For each equilibrium,  assuming that the $H^2$ norm  of  the initial data  $(f_0,B_0,U_0)$ is sufficient small  and  $f_0(x,v)$  has a compact support in the position $x$
and the velocity $v$, we construct the global well-posedness and decay rates of 
strong solutions near the equilibrium  in the whole space $\mathbb{R}^3$.
And the solution decays   polynomially. The global existence  result still holds for the  torus $\mathbb{T}^3$ case without the compact support assumption in   $x$. In addition, the decay rates are  
exponential. Lack  of dissipation structure in the Vlasov equation and the strong 
trilinear coupling term $((u-v)\times B)f$ in the model  are two main impediments in obtaining  our results. 
To surround these difficulties, we assume   that $f_0(x,v)$  has a compact support and 
  utilize the method of characteristics to calculate the size of the supports of $f$. Thus, we overcome 
the difficulty in estimating the integration   $\int_{\mathbb{R}^3}
\big((u-v)\times B\big)f\mathrm{d}v$  and obtain the global existence of strong solutions
  by taking advantage of a refined energy method. Moreover,  by making full use of the Fourier techniques,
we obtain the optimal time decay rate of the gradient of the solutions. This is the first
result on strong solutions to the Vlasov-MHD model containing nonlinear Lorentz forces.
\end{abstract}



\maketitle

\setcounter{equation}{0}
 \indent \allowdisplaybreaks

\section{Introduction and main result}\label{Sec:intro-resul}

 \subsection{Previous literature and our model}


Fluid-particle interaction models usually consist of a kinetic equation for a cloud of 
dispersed suspended particles in a fluid
and a hydrodynamic system for the  fluid. They 
have been widely used in varieties fields such as  sedimentation
of solid particles by external
force\cite{BD-jhde-2006},
combustion theory\cite{WF-th-1985}, chemical engineering \cite{BBJ-esaimp-2005}, rain formation \cite{FFS-02},
and medicine treatment\cite{Ma-thesis-2009}.

Generally speaking, there are many  factors affecting  the modelling of the fluid-particle interaction process. These factors include, for example,  the incompressibility/compressibility of the fluid, 
 the temperature effect  of  the fluid, 
the properties of the particles, the domain (bounded, periodic, or the whole space) occupied by the fluid, the interaction form of the particles and the fluid, and the effect of Brownian of particles.  A common fluid-particle coupled model is the following 
incompressible Vlasov-Navier-Stokes system: 
\begin{equation}\label{I-1''}
\left\{
\begin{aligned}
& \partial_t F+v \cdot \nabla_x F-{\rm div}_v((v-u) \cdot F)=0, \\
& \partial_t u+u\cdot \nabla_x u+\nabla_x P-\mu \Delta_x u=\int_{\mathbb{R}^3} 
(v-u)F\,\mathrm{d}v, \\
& \nabla_x \cdot u=0.
\end{aligned}\right.
\end{equation}
Here the unknown $f(t, x, v)$ denotes the distribution function of particles
at time $t \geq 0$, with the position $x\in \Omega \subseteq \mathbb{R}^3$ 
and the velocity $v=\left(v_1, v_2, v_3\right) \in \mathbb{R}^3$.
 And  $u(t,x)$  denotes  a velocity field of the fluid and  
  $P(t,x)$  the pressure of the fluid. The particles and the fluid interact mutually via the so-called Brinkman force.  
 Boudin et al. \cite{BDGM-die-2009} established the global weak solution
to the system \eqref{I-1''}  on a   torus $\mathbb{T}^3$. Yu \cite{Yc-jmpa-2013} obtained the global existence
of weak solutions to \eqref{I-1''}  in a bounded domain $\Omega\subset \mathbb{R}^3$. 
Choi and Kwon \cite{CK-nonlinearity-2015} obtained the global well-posedness and large-time
behavior for the inhomogeneous Vlasov–Navier–Stokes equations on $\mathbb{T}^3$, where 
the solution exists on $[0,T]$ for    arbitrarily large $T>0$.  
Han-Kwan et al. \cite{HMM-arma-2020, EHM-N-2021, Hd-pmp-2022} obtained  the global existence and large time behavior of
weak solutions to \eqref{I-1''}  near the trivial equilibrium $(\bar f,\bar u)=(0,0)$ 
on a torus $\mathbb{T}^3$, the whole space $\mathbb{R}^3$, and a  bounded domain $\Omega \subset\mathbb{R}^3$, respectively. Very recently, Danchin \cite{Danchin-24} obtained the   global  strong solution of  Fujita-Kato type  to \eqref{I-1''}  in the whole space $\mathbb{R}^3$ and the  optimal rate of it.
It should be pointed out that although there is no apparent dissipation structure in the Vlasov equation \eqref{I-1''}$_1$, the system \eqref{I-1''} definitely has a dissipation effect  term $\frac{1}{2}\iint_{\mathbb{R}^3_x\times \mathbb{R}^3_v}|v -u|^2f{\rm d}v{\rm d}x$ in its energy equality:
\begin{align*}
 \frac12\frac{\rm d}{{\rm d}t} \Big\{\iint_{\mathbb{R}^3_v\times \mathbb{R}^3_x}
|v|^2f{\rm d}v{\rm d}x+   \int_{\mathbb{R}_x^3}| u|^2{\rm d}x\Big\}
+\int_{\mathbb{R}^3_x}
\mu |\nabla u|^2{\rm d}x+ \iint_{\mathbb{R}_v^3\times \mathbb{R}_x^3}
|v -u|^2f{\rm d}v {\rm d}x=0, 
\end{align*}
  which plays a crucial role in the arguments of  
the decay rate of solutions \cite{Hd-pmp-2022, Danchin-24}.
When the Brownian motion effect is considered, i.e. the equation \eqref{I-1''}$_1$ includes  a diffusion term 
$\nu \Delta_v f$ and reads as 
\begin{align}\label{fpa}
 \partial_t F+v \cdot \nabla_x F=-\mathcal{L}F,  
\end{align}
where $  \mathcal{L}F$ is a    Fokker-Planck type operator defined as $\mathcal{L}F=\nabla_v \cdot\left(\nabla_v F+(v-u) F\right)$, 
then the global Maxwellian  
\begin{align}\label{maa}
  M=M(v)=(2\pi)^{- \frac{3}{2}}e^{- \frac{|v|^2}{2}}, 
\end{align}  
together with $\bar u=0$, is an equilibrium  to the new system  \eqref{fpa}, \eqref{I-1''}$_2$,  and \eqref{I-1''}$_3$. In this situation, it is much easier to obtain 
the global classical solution to it, see \cite{GHMZ-sjma-2010} on the torus case   and  \cite{CKL-jhde-2013, DL-krm-2013, LMW-sjma-2017, MW-cvpde-2020} for more complicated compressible  fluid-particle models.

In plasma physics, to describe the interaction of the magnetohydrodynamic (MHD) bulk fluid with a rarified 
ensemble of energetic particles, many kinetic-MHD models have been presented and some numerical simulations are given  \cite{cheng-91,HT-2011, T-2010,TT-2015}. 
Based on the so-called current-coupling scheme, Cheng, S\"{u}li and Tronci \cite{CST-plms-2017}  introduced the 
 following   Vlasov-MHD model: 
\begin{equation}\label{I1}
\left\{
\begin{aligned}
& \partial_t f+v \cdot \nabla_x f+\big((v-u)\times B\big) \cdot \nabla_v f=0, \\
& \partial_t u+u\cdot \nabla_x u+\nabla_x P-(\nabla_x \times B)\times B-\mu_1\Delta_x u=\int_{\mathbb{R}^3} 
\big((u-v)\times B\big)f\mathrm{d}v, \\
& \partial_t B-\nabla_x \times(u\times B)-\mu_2\Delta_x B=0, \\
& \nabla_x \cdot u=0, \quad \nabla_x \cdot B=0,
\end{aligned}\right.
\end{equation}
with the initial data
\begin{equation}\label{I2}
f(0, x, v)=f_0(x, v), \quad u(0, x)=u_0(x), \quad B(0, x)=B_0(x) 
\end{equation}
on a torus $\mathbb{T}^3$. As before, the unkown $f(t, x, v)$ denotes the distribution function of particles
at time $t \geq 0$, with the position $x=\left(x_1, x_2, x_3\right)\in  \mathbb{T}^3$ 
and the velocity $v=\left(v_1, v_2, v_3\right) \in \mathbb{R}^3$.
$u(t,x)$ denotes the velocity field of the fluid and $B(t,x)$   the magnetic field. 
And  $P(t,x)$ is the pressure of the fluid. The constant $\mu_1>0$ is the fluid viscosity coefficient,
and
$\mu_2>0$ is the magnetic diffusion coefficient. The particle distribution function $f$ couples with the magnetic field $B$ and the  velocity field $u$ via  the  nonlinear Lorentz force  $(v-u)\times B$.
Cheng, S\"{u}li and Tronci  \cite{CST-plms-2017}  established the global existence of finite energy weak solutions 
to the    Vlasov-MHD model  \eqref{I1} on $\mathbb{T}^3$ with  large initial data.
  
We point out that the fluid-particle model \eqref{I-1''} can be not derived formally by ignoring the magnetic field
 $B$ in the system \eqref{I1}. In fact, if we take $B\equiv 0$ in \eqref{I1}, then the terms involving the Brinkman force disappear and 
 the system is decoupled. The system  \eqref{I1}  has an energy conservation identity 
 \begin{align*}
 \frac12\frac{\rm d}{{\rm d}t} \Big\{\iint_{\mathbb{R}^3_v\times \mathbb{R}^3_x}
|v|^2f{\rm d}v{\rm d}x+   \int_{\mathbb{R}^3_x}(| u|^2+|B|^2){\rm d}x\Big\}
+\int_{\mathbb{R}^3_x}
\mu_1 |\nabla u|^2+\mu_2 |\nabla B|^2 {\rm d}x =0. 
\end{align*}
 In stark contrast to the fluid-particle model \eqref{I-1''}, there is no dissipation  effect 
caused by   $\frac{1}{2}\iint_{\mathbb{R}^3_x\times \mathbb{R}^3_v}|v -u|^2f{\rm d}v{\rm d}x$, which 
brings some difficulties  in proving  the global existence of strong solutions. 
 
 We mention that there are some results on   variety versions of the incompressible Vlasov-MHD model 
 \eqref{I1}. Chen, Hu and Wang \cite{CHW-jmfm-2016} obtained the global existence of weak solutions to the 
incompressible inhomogeneous Vlasov-MHD equations in a 3D bounded domain where
the particles and fluid coupling term doesn't involving the effect of magnetic field,
  i.e. they took the drag force $k\rho(v-u)$ instead of  the Lorentz force $ (v-u)\times B$ in the integral term 
  of \eqref{I1}$_2$.   Later,  Jiang \cite{Jp-jde-2017} obtained the global existence and time decay of smooth solutions  to the compressible Vlasov-Fokker-Planck-MHD system 
 in the whole space $\mathbb{R}^3$ with  the drag force 
 {which is supposed to be proportional to 
  the relative velocity $ (v-u)$}
  in the integral term.  Recently,    the first two authors of this paper \cite{LN-24} established 
  the global existence and time decay of smooth solutions to the 
  system \eqref{I1} by adding a  Fokker-Planck operator $\nabla_v \cdot (\nabla_v F+v F )$.
 
 To our best knowledge, there is no result on strong solutions to  the system \eqref{I1}.  
 It is readily to verify that  \eqref{I1} has two equilibria: 
 \begin{align}
  (\bar f,\bar u,\bar B)=(0,0,0)  \quad \text{and} \quad ( \tilde f,\tilde u,\tilde B)=(M,0,0), \label{teq}
 \end{align}
 where $M$ is the global Maxwellian defined in \eqref{maa}.   In this paper, we shall
  study the global well-posedness and decay rates of 
strong solutions to  \eqref {I1} near the foregoing two  equilibria in the whole space $\mathbb{R}^3$ 
and the torus $\mathbb{T}^3$, respectively.

\subsection{Main results}

In this subsection, we present our results  to the Cauchy problem \eqref{I1}--\eqref{I2} when the initial data
is a small perturbation
around the  aforementioned   two spatially steady states   $(\bar f,\bar u,\bar B) $ and $( \tilde f,\tilde u,\tilde B) $. We use $\Omega$ to denote  $\mathbb{R}^3_x$ 
or  $\mathbb{T}^3_x$. Due to the fact that the coefficients $\mu_1$ and $\mu_2$ 
have no effect on our results, 
 we  
set $\mu_1=\mu_2=1$ for presentation simplicity. Applying   the following identities:
\begin{align}
(\nabla_x \times B)\times B&=B\cdot\nabla_x B-\frac{1}{2}\nabla_x(|B|^2),\label{N1.4}\\
\nabla_x \times(u\times B)&=u\nabla_x\cdot B-u\cdot\nabla_x B+B\cdot \nabla_x u-B\nabla_x\cdot u, \label{N1.5}
\end{align}
 we  rewrite the system \eqref{I1} as follows:
\begin{equation}\label{I3}\left\{
\begin{aligned}
& \partial_t f+v \cdot \nabla_x f+\big((v-u) \times B\big) \cdot \nabla_v f =0, \\
& \partial_t u+u\cdot \nabla_x u+\nabla_x P-\Delta_x u-B\cdot\nabla_x B+\frac{1}{2}\nabla_x (|B|^2)=\int_{\mathbb{R}^3} 
\big((u-v)\times B\big)f\mathrm{d}v, \\
& \partial_t B-\Delta_x B+u\cdot\nabla_x B-B\cdot\nabla_x u=0, \\
& \nabla_x \cdot u=0, \quad \nabla_x \cdot B=0,
\end{aligned}\right.
\end{equation}
with the initial data 
\begin{align}\label{I--3}
(f(t, x, v),   u(t, x), B(t, x))|_{t=0}=(f_0(x, v),  u_0(x),  B_0(x)),\quad   x\in \Omega.
\end{align}

Before stating our results, we introduce some notions. 
The symbol $C$ represents a generic constant that may vary from line to line,
$A\sim B$ means that $\frac{A}{C}\leq B\leq CA$ for some positive constant $C$.
For an integrable function $g: \mathbb{R}_x^3 \rightarrow \mathbb{R}$, we define the Fourier transform as follows
\begin{align*}
\hat{g}(\xi)=\mathcal{F} g(\xi)=\int_{\mathbb{R}^3} \mathrm{e}^{-2 \pi i x \cdot \xi} g(x) \,\mathrm{d} x,
\end{align*}
where $x \cdot \xi=\sum_{j=1}^3 x_j \xi_j$ for $\xi \in \mathbb{R}^3$, and $i=\sqrt{-1} \in \mathbb{C}$ is the imaginary unit. And $g \cdot \bar{h}=(g \mid h)$ denotes the dot product of $f$ with the complex conjugate of $ h $.

For brevity, we use $\langle\cdot, \cdot\rangle$ to denote  the standard $L^2$ inner product in $\mathbb{R}_v^3$, and $(\cdot, \cdot)$   the $L^2$ inner product in $\Omega$ or $\Omega \times \mathbb{R}_v^3$, and the corresponding norms can be defined by $|\cdot|_2$ and $\|\cdot\|$, respectively.
For   the multi-indices $\alpha=\left(\alpha_1, \alpha_2, \alpha_3\right)$ and $\beta=\left(\beta_1, \beta_2, \beta_3\right)$,  we denote
$$
\partial_\beta^\alpha=\partial_{x_1}^{\alpha_1} \partial_{x_2}^{\alpha_2}
\partial_{x_3}^{\alpha_3} \partial_{v_1}^{\beta_1} \partial_{v_2}^{\beta_2} \partial_{v_3}^{\beta_3}, 
$$
and the norms 
\begin{align*}
\|w\|_{H^s}:=\sum_{|\alpha| \leq s}\|\partial^\alpha  w\|, \quad\|g\|_{H_{x, v}^s}:=\sum_{|\alpha|+|\beta| \leq s}\|\partial_\beta^\alpha g\|.
\end{align*}
Finally, we use  $C_\alpha^\beta$ to denote the general binomial coefficient. 
 
\subsubsection{The case $(\bar f,\bar u,\bar B)=(0,0,0)$} For this equilibrium  state, we have 
 
\begin{thm}\label{I4}
Let $\Omega=\mathbb{R}^3$ and suppose that the initial data $(f_0,u_0,B_0)$ satisfy 
$f_0\in H_{x,v}^{2}$, $u_0,B_0\in H^2$,
and  $f_0(x,v)$ has a compact 
support in $x$ and $v$.
Then, for  any time $T\in(0,+\infty)$,   
  there
exist two constants $\epsilon_1=\epsilon_1(T)>0$ and $\epsilon_2=\epsilon_2(T)>0$, such that
if $\|f_0\|_{H_{x,v}^2}<\epsilon_1$ and 
$\|(u_0,B_0)\|_{H^{2}}<\epsilon_2$, then the Cauchy problem \eqref{I3}--\eqref{I--3} admits a unique global strong solution $(f, u,B)$ on $[0,T]$ satisfying 
\begin{gather}
f\in C \big([0,T];H_{x,v}^{2}\big),\    u,B\in C\big([0,T];H^{2}\big)\cap L^2\big(0,T;H^{3}\big);\label{t3333}\\ 
 \sup_{0\leq t \leq T}\|(u,B)(t)\|_{H^{2}}\leq C\|(u_0,B_0)\|_{H^{2}};\label{KKK}\\
 \sup_{0\leq t \leq T}\|f(t)\|_{H_{x,v}^2} {\leq C \|f_0\|_{H_{x,v}^2}^{\frac{1}{2}}}.\label{G1.5}
\end{gather}
Moreover, if we further assume that
$\|(u_0,B_0)\|_{L_1}<\varsigma$ with $\varsigma \in (0,1)$, then
\begin{align}
\|(u,B)(t)\|_{{H}^{2}}& \leq C(1+t)^{-1}, \label{G1.6}\\
\|\nabla(u,B)(t)\|_{{H}^{1}}& \leq C(1+t)^{-\frac{5}{4}}  \label{G1.7}
\end{align}
for any $t\in[0,T]$.
\end{thm}

\begin{rem} \label{reaa}
 In Theorem \ref{I4}, we can choose  $0<\epsilon_1\ll 1$  with an order of $e^{-\mathcal{O}(1)T}$,
and $0<\epsilon_2\ll 1$ with an order of $\frac{\mathcal{O}(1)}{T+1}$. 
\end{rem}

\begin{rem}
Due to the initial condition
$\|(u_0,B_0)\|_{H^{2}}<\epsilon_1$, we 
  obtain the time decay rate  of our solutions with $(1+t)^{-1}$ which is faster than the usual  optimal time decay rate $(1+t)^{-3/4}$ for the incompressible MHD equations \cite{SSS-1996} in the whole space $\mathbb{R}^3$. 
Besides, we  also obtain the optimal time decay 
rates with regard to gradient of the solutions $u$ and $B$ in the sense of that 
they coincide with that of the  linearized MHD equations. 
\end{rem}

\begin{thm}\label{I5}
Let $\Omega=\mathbb{T}^3$  and suppose that the initial data $(f_0,u_0,B_0)$ satisfy 
 $f_0\in H_{x,v}^{2}$, $u_0,B_0\in H^2$,  $f_0(x,v)$ has a compact 
support in   $v$, and  $B_0 $ satisfies 
\begin{align*}
& \int_{\mathbb{T}^3}B_0(x)\mathrm{d}x=0.
\end{align*} 
Then, for  any time $T\in(0,+\infty)$,   
  there exist $\epsilon_3=\epsilon_3(T)>0$ and $\epsilon_4=\epsilon_4(T)>0$, such that if 
$\|f_0\|_{H_{x,v}^2}<\epsilon_3$ and
$\|(u_0,B_0)\|_{H^{2}}<\epsilon_4$, then
  the Cauchy problem \eqref{I3}--\eqref{I--3} admits a unique strong solution $(f, u, B )$ on $[0,T]$ 
   satisfying 
\begin{gather}
  f\in C \big([0,T];H_{x,v}^{2}\big) ,\ u,B\in C\big([0,T];H^{2}\big)\cap L^2\big(0,T;H^{3}\big); \nonumber\\
  \|( u, B)(t)\|_{H^2} \leq C e^{-\lambda_a t} \label{peri-a} 
\end{gather}
for any  $t\in[0,T]$, where $\lambda_a>0$ is a constant independent of $t$.
\end{thm}


\begin{rem}
As in Remark \ref{reaa}, we can choose  $0<\epsilon_3\ll 1$ with  an order of $e^{-\mathcal{O}(1)T}$,
and $0<\epsilon_4\ll 1$ with  an order of $\frac{\mathcal{O}(1)}{T+1}$.
\end{rem}

\begin{rem}
Considering the boundedness of the periodic domain, 
here we only need to  assume that $f_0(x,v)$ has a compact 
support  in the velocity $v$.  
\end{rem}

\begin{rem} \label{reaab}
Theorems \ref{I4} and \ref{I5} still hold for the incompressible Vlasov–Navier–Stokes equations considered in \cite{CK-nonlinearity-2015},  where 
 the authors chose the initial data $f_0$ and $u_0$  
 sufficient small satisfying  the order of $e^{-\mathcal{O}(1)T}$.  
Here  we  improve their result by relaxing the assumption 
on the fluid velocity $u_0$ to an order of  $\frac{\mathcal{O}(1)}{T+1}$. Moreover, we obtain 
the time decay of the solutions and their gradients. 
\end{rem}

\subsubsection{The case $( \tilde f,\tilde u,\tilde B)=(M,0,0)$}

In this subsection, we shall prove
the global existence of strong solutions
to the Cauchy problem \eqref{I1}--\eqref{I2}
when the initial data are a small perturbation
near the  equilibrium state 
$(\tilde f,\tilde u,\tilde B)=(M,0,0)$.
Setting $f=g+M$,  
 we   rewrite the system \eqref{I1} through \eqref{N1.4}--\eqref{N1.5} as follows (taking $\mu_1=\mu_2=1$) :
\begin{equation}\label{I8}\left\{
\begin{aligned}
& \partial_t g+v \cdot \nabla_x g+\big((v-u) \times B\big) \cdot \nabla_v g +(u\times B)\cdot v M=0, \\
& \partial_t u+u\cdot \nabla_x u+\nabla_x P-\Delta_x u-B\cdot\nabla_x B\\
& \qquad  +\frac{1}{2}\nabla_x (|B|^2)-u\times B=\int_{\mathbb{R}^3} 
\big((u-v)\times B\big)g\mathrm{d}v, \\
& \partial_t B-\Delta_x B+u\cdot\nabla_x B-B\cdot\nabla_x u=0, \\
& \nabla_x \cdot u=0, \quad \nabla_x \cdot B=0,
\end{aligned}\right.
\end{equation}
with the initial data 
\begin{align}\label{I--8}
(g(0, x, v),   u(t, x), B(t, x))|_{t=0}=\,&((f_0(x, v)-M,  u_0(x),  B_0(x)) \nonumber\\
:=\,& (g_0(x, v),  u_0(x),  B_0(x)),   \quad x\in \Omega.
\end{align} 

We have the following results.

\begin{thm}\label{T5.2}
Let $\Omega=\mathbb{R}^3$ and suppose that the initial data $(g_0,u_0,B_0)$ satisfy  
$g_0\in H_{x,v}^{2}$, $u_0,B_0\in H^2$,
and  $g_0(x,v)$ has a compact 
support in $x$ and $v$.
Then, for any time $T\in(0,+\infty)$,   there
exists a constant $\tau_1=\tau_1(T)>0$  such that if 
$\|g_0\|_{H_{x,v}^2}+\|(u_0,B_0)\|_{H^{2}}<\tau_1$, then 
the Cauchy problem \eqref{I8}--\eqref{I--8} admits a unique global strong solution $(g, u,B)$ 
on $[0,T]$ satisfying 
\begin{gather}
 g\in C \big([0,T];H_{x,v}^{2}\big),\ u,B\in C\big([0,T];H^{2}\big)\cap L^2\big(0,T;H^{3}\big);\nonumber\\
 \sup_{0\leq t \leq T}\|(u,B)\|_{H^{2}}\leq C\|(u_0,B_0)\|_{H^{2}};\label{G51}\\
\label{G5.3}
 \sup_{0\leq t \leq T}\|g(t)\|_{H_{x,v}^2}\leq C\big\{\|g_0\|_{H_{x,v}^2}+\|(u_0,B_0)\|_{H^2}\big\}^{\frac{1}{2}}.
\end{gather}
\end{thm}

\begin{thm}\label{T5.1}
Let $\Omega=\mathbb{T}^3$ and suppose that the initial data $(g_0,u_0,B_0)$ satisfy 
 $g_0\in H_{x,v}^{2}$, $u_0,B_0\in H^2$,    $g_0(x,v)$ has a compact 
support in $v$, and  
\begin{gather*}
  \int_{\mathbb{T}^3\times\mathbb{R}^3} g_0(x,v) \mathrm{d} x\mathrm{d}v=0,  \ \ \int_{\mathbb{T}^3}B_0(x)\mathrm{d}x=0,\\
  \int_{\mathbb{T}^3\times\mathbb{R}^3} vg_0(x,v) \mathrm{d} x\mathrm{d}v+\int_{\mathbb{T}^3} u_0(x) \mathrm{d} x=0.
\end{gather*}
Then, for any $T\in(0,+\infty)$,   there
exists  a constant $\tau_2=\tau_2(T)>0$  such that if
$\|g_0\|_{H_{x,v}^2}+\|(u_0,B_0)\|_{H^{2}}<\tau_2$, 
Then, the Cauchy problem \eqref{I8}--\eqref{I--8} admits a unique strong solution $(g, u, B )$ 
on $[0,T]$ satisfying 
\begin{gather}\label{peri-b}
 g\in C \big([0,T];H_{x,v}^{2}\big),\    u,B\in C\big([0,T];H^{2}\big)\cap L^2\big(0,T;H^{3}\big); \nonumber \\
  \|( u, B)(t)\|_{H^2} \leq C e^{-\lambda_{b} t} 
\end{gather}
for any  $t\in[0,T]$, where $\lambda_b>0$ is a constant independent of $t$.
\end{thm}

\begin{rem}
In Theorems \ref{T5.2} and \ref{T5.1},  we can choose $0<\tau_1\ll 1$ and  $0<\tau_2\ll 1$ have  an order of $e^{-\mathcal{O}(1)T}$.
\end{rem}

\begin{rem}
Compared to the vacuum case   $(\bar f,\bar u,\bar B)=(0,0,0)$,
  for the case $( \tilde f,\tilde u,\tilde B)=(M,0,0)$,
   there are two  additional terms in the system \eqref{I8}, i.e.  $(u\times B)\cdot v M$ in the equation of $g$ 
    and $ u\times B$ in the equation of $u$, which bring  some difficulties in handling  the  upper bound of $\|g\|_{H^2_{x,v}}$ in the energy estimates, i.e.  
$C\|u\|_{H^2}\|B\|_{H^2}$, which may result in that $u$ or $B$ has an initial value controlled by an order of $e^{-\mathcal{O}(1)T}$. Combining the energy estimate part, we need to assume that the initial values of $u$ and $B$ are controlled by an order of $e^{-\mathcal{O}(1)T}$. 
Notice that exponential decay occurs faster than algebraic decay. Due to the algebraic decay obtained by performing Fourier transform on the whole space, we don't need to consider the large time behavior here.
\end{rem}

\begin{rem}
  It should be pointed out that for the fluid-particle model \eqref{I-1''},  the global Maxwellian $M$ 
  is not an  equilibrium with $\bar u=0$. This is one of main differences between the fluid-particle model and 
  the Vlasov-MHD model. 
  \end{rem}


Now we give some comments on the proofs of our results. In order to obtain the global existence of strong solutions, we first construct  a linearized system to achieve the existence and  uniqueness of solutions.  Then, we shall establish the existence of the strong solutions
of the original nonlinear system  through constructing the approximation sequences. 
 Compared to the Vlasov-Fokker-Planck-MHD system studied  in  \cite{LN-24} 
 where  the refined energy estimates were obtained  by applying  the micro-macro  
decomposition of particles,  for our model \eqref{I3}--\eqref{I--3},  this method is ineffective because of  the   lack  of dissipation structure containing in the Fokker-Planck operator $\mathcal{L}f$. In addition,  without the   strong dissipation structure, the 
trilinear coupling term $((u-v)\times B)f$  brings some extra troubles in  energy estimates. 
 To surround these difficulties, inspired by the ideas developed in \cite{CK-nonlinearity-2015},  we assume   that $f_0(x,v)$  has a compact support in $x$ and $v$ and 
  utilize the method
of characteristics to calculate the size of
the supports of $f$. 
By establishing a priori estimates for $u$ and $B$ respectively  (see \eqref{HH2.1}),
  we obtain a more precise estimate for the size of the supports of $f$ in the position and velocity
and overcome 
the difficulty in estimating the
integration   $\int_{\mathbb{R}^3}
\big((u-v)\times B\big)f\mathrm{d}v$. With these preparations, we   take advantage of
a refined energy method to enclose the a priori estimates and 
finally obtain the global existence of strong solutions. 
In these process, we only require that the initial values of $u,B$ can be controlled by an order of $\frac{\mathcal{O}(1)}{T+1}$ which  is different from \cite{CK-nonlinearity-2015} that the initial values of $u,B$ should be controlled by an order of $e^{-\mathcal{O}(1)T}$.
Moreover,  by making full use of the Fourier techniques,
we obtain the optimal time decay
rate of the gradient of the solutions. For the periodic case, due to lack of the  dissipation operator 
$\mathcal{L}f$ and the equation of $u$ involving $f$, we encounter  additional difficulties 
in obtaining the exponential decay of $u$ and $B$. 
We ultimately solve it by constructing a differential inequality  (see \eqref{K4.7} below).

When the perturbation is near the global Maxwellian $M$,  
the new term $u\times B$  appears in the velocity equation and  produces some troubles.
Fortunately, for the ``high-order”  estimates involving  this term, we   utilize the energy method to handle them. 
For  the ``zero-order" estimate concerning it, due to the fact that we can't shrink the term $\int_{\mathbb{R}^3}u\cdot(u\times B)\mathrm{d}x$ into the form with energy multiplied by the square of dissipation, we finally  use the identity $u\cdot(u\times B)=-B\cdot(u\times u)=0$  to
overcome it. 
In addition, the new term $(u\times B)\cdot v M$  
also bring some difficulties for the estimate of $F$.
To solve these difficulties,  we   limit some conditions on the initial data,
which implies that the both $u$ and $B$ have exponential decay
rate. To the best of our knowledge, the  results obtained in this paper are the 
first ones  on strong solutions to the Vlasov-MHD model containing nonlinear Lorentz forces.

The rest of this paper is organized as follows.  In Section 2, we first establish the 
existence and uniqueness of strong solutions for the
linearized problem \eqref{H2.1}--\eqref{I--33} of the nonlinear problem \eqref{I3}--\eqref{I--3} near 
the equilibrium  $(\bar f,\bar u,\bar B)=(0,0,0)$ via the method of characteristics and  the  energy method. By doing this,
we  establish the global existence and uniqueness of 
strong solutions to  \eqref{I3}--\eqref{I--3}. We also obtain the upper bounds   of strong solutions.
In Section 3, we devote ourselves to the proofs of \eqref{G1.6}--\eqref{G1.7} in Theorem \ref{I4}. Using
the  upper bounds  of strong solutions and 
 developing decay estimates on the linearized equations  of \eqref{I3}, 
we further obtain the optimal decay rates of $u, B$ and
the gradient of them in $\mathbb{R}^3$.
In Section 4, we discuss the strong solutions to the system \eqref{I3}--\eqref{I--3} in 
a torus $\mathbb{T}^3$ and 
achieve the exponential decay of $u$ and $B$ near the equilibrium $(\bar f,\bar u,\bar B)=(0,0,0)$. 
Finally, in Section 5,
 we obtain the global well-posedness and decay rates of strong solutions to 
 the   problem \eqref{I3}--\eqref{I--3} near the 
 equilibrium $( \tilde f,\tilde u,\tilde B)=(M,0,0)$ in both  $\mathbb{R}^3$ and $\mathbb{T}^3$.

\section{Global existence of strong solutions   to the problem \eqref{I3}--\eqref{I--3} in  $\mathbb{R}^3$} 

In this section, we first establish the existence 
of strong solutions to the linearized problem  of the 
nonlinear   problem \eqref{I3}--\eqref{I--3} near  
the equilibrium  state $(\bar f,\bar u,\bar B)=(0,0,0)$ in $\mathbb{R}^3$. Then we obtain   the existence and uniqueness of strong solutions to the nonlinear problem.

\smallskip 
To begin with,  we recall  some basic facts which   will be used frequently.  

\begin{lem}\label{L1}
There exist positive constants C, such as for any $f,g \in H^2\left(\mathbb{R}^3\right)$ and any multi-index $\alpha$ verifying $1\leq |\alpha|\leq 2$, one has
\begin{align*}
\|f\|_{L^{\infty}\left(\mathbb{R}^3\right)} 
\leq&\, C\|\nabla_x f\|_{L^{2}\left(\mathbb{R}^3\right)}^{\frac{1}{2}}\left\|\nabla_x^2 f\right\|_{L^{2}\left(\mathbb{R}^3\right)}^{\frac{1}{2}},\\ \|fg\|_{H^{2}\left(\mathbb{R}^3\right)}
\leq& \, C\|f\|_{H^{2}\left(\mathbb{R}^3\right)}\|\nabla_x g\|_{H^{2}\left(\mathbb{R}^3\right)},\\
\|\partial^{\alpha}fg\|_{L^2\left(\mathbb{R}^3\right)}
\leq & \,  C\|\nabla_x f\|_{H^{1}\left(\mathbb{R}^3\right)}\|\nabla_x g\|_{H^{1}\left(\mathbb{R}^3\right)} .
\end{align*}
\end{lem}
The proof  of this    lemma   can be   found in
\cite{AF-Pa-2003, Dk-mz-1992} 
(see Lemmas 2.1 and 2.2 in \cite{Dk-mz-1992} and Lemma 2.1
in \cite{Dk-mz-1992}).

\begin{lem}[\!\cite{CDM-krm-2011}]\label{H2}
Given any $0<\beta_1\neq1$
and $\beta_2>1$, it holds that
\begin{align*}
\int_0^t(1+t-s)^{-\beta_1}(1+s)^{-\beta_2}\mathrm{d}s\leq C(1+t)^{-\min\{\beta_1,\beta_2\}} 
\end{align*}
for all $t\geq 0$.
\end{lem}

\subsection{Global existence and uniqueness of a linearized system} 
Now we prove the existence and uniqueness of the strong solution
for the linearized system as follows:
\begin{equation}\label{H2.1}\left\{
\begin{aligned} 
& \partial_t f+v \cdot \nabla_x f+\big((v- \mathsf{u}) \times \mathsf{B}\big) \cdot \nabla_v f =0, \\
& \partial_t u+\mathsf{u}    \cdot \nabla_x u+\nabla_x P-\Delta_x u-\mathsf{B}\cdot\nabla_x B+\frac{1}{2}\sum_{i,j=1}^3 \mathsf{B}_i  \partial_j B_i, \\
&\qquad  =\int_{\mathbb{R}^3} (u\times \mathsf{B})f\mathrm{d}v-\int_{\mathbb{R}^3} (v\times B)f\mathrm{d}v,\\
& \partial_t B-\Delta_x B+\mathsf{u} \cdot\nabla_x B-\mathsf{B}\cdot\nabla_x u=0, \\
& \nabla_x \cdot u=0, \quad \nabla_x \cdot B=0 
\end{aligned}\right.
\end{equation}
with the initial data 
\begin{align}\label{I--33}
(f(t, x, v),   u(t, x), B(t, x))|_{t=0}=(f_0(x, v),  u_0(x),  B_0(x)),\quad   x\in \mathbb{R}^3,
\end{align} 
where $\mathsf{u}$ is a known velocity field and $\mathsf{B}$ is
a known magnetic field.

\begin{thm}\label{reaaa}
Suppose that  $f_0(x,v)$ has a compact 
support in $x$ and $v$.
Then, for  any time $T\in(0,+\infty)$,   
  there
exist two constants $\epsilon_1=\epsilon_1(T)>0$ and $\epsilon_2=\epsilon_2(T)>0$, such that
if $\|f_0\|_{H_{x,v}^2}<\epsilon_1$ and
$\|(u_0,B_0)\|_{H^{2}}<\epsilon_2$, and $\mathsf{u}$, $\mathsf{B}$ satisfy the following smallness and regularity conditions:
\begin{align*}
 \|(\mathsf{u} ,\mathsf{B})\|_{C ([0,T];H^{2})}<\epsilon_2,\quad 
 \|(\mathsf{u} ,\mathsf{B})\|_{L^2 (0,T;H^{3})}<\epsilon_2.
\end{align*}
Then the problem \eqref{H2.1}--\eqref{I--33} enjoys a unique strong solution on $[0,T]$ such that 
\begin{align}\label{t2333}
\|f\|_{C ([0,T];H_{x,v}^{2})}<\epsilon_1^{\frac{1}{3}}, \quad\|(u,B)\|_{C ([0,T];H^{2})}<\epsilon_2,\quad 
\|(u,B)\|_{L^2 (0,T;H^{3})}<\epsilon_2.
\end{align}
\end{thm}

\begin{rem} 
 In Theorem \ref{reaaa}, we can choose  $0<\epsilon_1\ll 1$  with an order of $e^{-\mathcal{O}(1)T}$,
and $0<\epsilon_2\ll 1$ with an order of $\frac{\mathcal{O}(1)}{T+1}$. 
\end{rem}

To prove Theorem \ref{reaaa}, we divide it into several parts.

\subsubsection{The estimate of $f$}
Motivated by \cite{CK-nonlinearity-2015},
we solve the Vlasov-type equation
\eqref{I3}$_1$ for $f$ by the method of characteristics. For this, 
we need to estimate the size of the support
of $f$ in $x$ and $v$. 
Let $\Sigma_x(t)$ and $\Sigma_v(t)$ be the $x$-- and $v$--projections of supp$f(t,\cdot,\cdot)$ respectively. 
Below we define
\begin{align*}
\Sigma_x(t):=\,&\{x\in\mathbb{R}^3_x:\:\exists \, (x,v)\in\mathbb{R}^3_x\times\mathbb{R}^3_v \  \text{such that}\  f(t, x,v)\neq0\},\\
\Sigma_v(t):=\,&\{v\in\mathbb{R}^3_v:\:\exists \, (x,v)\in\mathbb{R}^3_x\times\mathbb{R}^3_v \  \text{such that}\  f(t, x,v)\neq0\}.
\end{align*}
Furthermore,
we  define $R_x( t) $ and $R_v( t)$  as follows:
\begin{align*}
\mathcal{R}_{x}(t):=\sup_{x\in\Sigma_{x}(t)}|x|,\qquad \mathcal{R}_{v}(t):=\sup_{v\in\Sigma_{v}(t)}|v| .
\end{align*}
For a given $(x,v)\in \mathbb{R}^3 \times\mathbb{R}^3$, 
we   define a forward trajectories $(X(s;0,x,v),V(s;0,x,v))$ as follows:
\begin{align}\label{G2.2}
\frac{\mathrm{d}X(s)}{\mathrm{d}s}&=V(s),\\ \label{G2.3}
\frac{\mathrm{d}V(s)}{\mathrm{d}s}
&= (V(s)-\mathsf{u} (s,X(s)) )\times \mathsf{B} (s,X(s)),\end{align}
with the initial condition satisfies $(X(0),V(0))=(x,v)$.
Below we give the growth estimates of the support of $f$ in $x$ and $v$.

\begin{lem}\label{N2.3}
For any $T>0$, suppose that $f_0(x,v)$ has a compact 
support in $x$ and $v$, then it holds  
\begin{align}
\mathcal{R}_{x}(t)&\leq \mathcal{R}_{x}(0)+C\left(\mathcal{R}_{v}(0)+\epsilon_2 T\right)T,\label{RX}\\
\mathcal{R}_{v}(t)&\leq C\big(\mathcal{R}_{v}(0)+\epsilon_2 T\big)
 \label{RV}
\end{align}
for all $0\leq t\leq T$.
\end{lem}

\begin{proof}
From \eqref{G2.3}, one   easily get
\begin{align*}
\frac{\mathrm{d}|V(s)|}{\mathrm{d}s}&\leq |V(s)|\|\mathsf{B} \|_{L^{\infty}}+\|\mathsf{u}\|_{L^{\infty}}\|\mathsf{B} \|_{L^{\infty}}\\
&\leq\|\mathsf{B} \|_{H^2}|V(s)|+\|\mathsf{u}\|_{H^2}\|\mathsf{B} \|_{H^2}\\
&\leq\epsilon_2|V(s)|+\epsilon_2^2\\
&\leq\epsilon_2|V(s)|+\epsilon_2.
\end{align*}
Here we have used the fact that $\epsilon_2<1$ is sufficient small.
By Gronwall’s inequality, we arrive at
\begin{align}\label{VTT}
|V(t)|&\leq e^{\epsilon_2 T}\big(|V(0)|+\epsilon_2 T\big) 
 \leq C\big(|V(0)|+\epsilon_2 T\big),
\end{align}
and hence we obtain \eqref{RV}. The estimate of $R_x(t)$ 
  follows from \eqref{VTT} and  
  \eqref{G2.2}.
\end{proof}
For the simplicity of symbols, we set
\begin{align*}
\mathcal{R}_{x}^{\infty}\,&:=\sup_{0\leq t\leq T}\mathcal{R}_{x}(t),\quad \mathcal{B}_x:=\mathcal{B}_{\mathcal{R}_{x}^{\infty}},
\\ \mathcal{R}_{v}^{\infty}\,&:=\sup_{0\leq t\leq T}\mathcal{R}_{v}(t),\quad \mathcal{B}_v:=\mathcal{B}_{\mathcal{R}_{v}^{\infty}},
\end{align*}
where $\mathcal{B}_{\mathcal{R}_{x}^{\infty}}$ and $\mathcal{B}_{\mathcal{R}_{v}^{\infty}}$ denote the balls with the  radii $\mathcal{R}_{x}^{\infty}$ and $\mathcal{R}_{v}^{\infty}$, respectively.

Thanks to the smallness of $\epsilon_1$, we obtain 
\begin{align}
\mathcal{R}_x^\infty&\approx C(1+T),\quad \ \mathcal{R}_v^\infty \approx C. \label{H2.6}
\end{align}

With the above preparation, we shall establish the estimate of $f$. 

\begin{lem}\label{L2.0}
For the chosen $T$ stated in Lemma \ref{N2.3}, the following inequality   
\begin{align}\label{G2.7}
&\frac{1}{2} \frac{\mathrm{d}}{\mathrm{d} t}\|f\|_{H_{x,v}^2}^2
\leq C\big(1+\|\nabla_x \mathsf{u}\|_{H^2}^2+\|\nabla_x \mathsf{B} \|_{H^2}^2\big)\|f\|_{H_{x,v}^2}^2
\end{align}
holds for all $0\leq t\leq T$.
\end{lem}
\begin{proof}
Firstly, we consider the zeroth-order estimate of $f$. 
Multiplying \eqref{H2.1}$_1$ by $f$, integrating the result over $\mathbb{R}^3 \times\mathbb{R}^3 $,
and noticing the fact that $ ((v-\mathsf{u}) \times \mathsf{B}  ) \cdot \nabla_v f
= {\rm div}_v  [((v-\mathsf{u}) \times \mathsf{B} )  f ]$,  one gets
\begin{align*}
\frac{1}{2} \frac{\mathrm{d}}{\mathrm{d} t}\|f\|^2=0 .
\end{align*}

Next, we study the pure spatial derivatives of $f$. 
Applying $\partial^{\alpha}$ with $1 \leq|\alpha| \leq 2$ 
to the equation \eqref{H2.1}$_1$, we have
\begin{align}\label{G2.8}
&\partial^\alpha \partial_tf+v\cdot\nabla_x\partial^\alpha f
+\sum_{|\alpha^\prime|\leq |\alpha|}C_{\alpha}^{\alpha^\prime}v\times\partial^{\alpha^\prime}\mathsf{B} \cdot\nabla_v
\partial^{\alpha-\alpha^\prime}f\nonumber\\
&  \qquad  -\sum_{|\alpha^\prime|\leq |\alpha|}C_{\alpha}^{\alpha^\prime}\partial^{\alpha^\prime}(\mathsf{u}\times \mathsf{B} )\cdot\nabla_v\partial^{\alpha-\alpha^\prime}f=0.  
\end{align}
Multiplying \eqref{G2.8}
by $\partial^{\alpha}f$
and taking integration over $\mathbb{R}^3\times\mathbb{R}^3$ gives
\begin{align*}
\frac{1}{2}\frac{\mathrm{d}}{\mathrm{d}t}\|\partial^\alpha f\|^2
= \,&\sum_{1\leq|\alpha^\prime|\leq |\alpha|}C_{\alpha}^{\alpha^\prime}
\int_{\mathbb{R}^3\times\mathbb{R}^3}
\big(v\times\partial^{\alpha^\prime}\mathsf{B} \cdot\nabla_v\partial^{\alpha-\alpha^\prime}f \big)\partial^\alpha f\mathrm{d}x\mathrm{d}v \\
&-\sum_{1\leq|\alpha^\prime|\leq |\alpha|}C_{\alpha}^{\alpha^\prime}\int_{\mathbb{R}^3\times\mathbb{R}^3}
\big(\partial^{\alpha^\prime}(\mathsf{u}\times \mathsf{B} )\cdot\nabla_v\partial^{\alpha-\alpha^\prime}f \big)\partial^\alpha f\mathrm{d}x\mathrm{d}v\\
\leq\,& C\int_{\mathbb{R}^3\times \mathsf{B} _{v}}\big(1+|\partial^{\alpha^\prime}\mathsf{B} |^2\big)|\nabla_v\partial^{\alpha-\alpha^\prime}f||\partial^\alpha f|\mathrm{d}x\mathrm{d}v\\
&+C\int_{\mathbb{R}^3\times \mathbb{R}^3}|\partial^{\alpha^\prime}(\mathsf{u}\times \mathsf{B} )||\nabla_v\partial^{\alpha-\alpha^\prime}f||\partial^\alpha f|\mathrm{d}x\mathrm{d}v\\
\leq\,& C\big(\|\nabla_x \mathsf{B} \|_{H^2}\| \mathsf{u}\|_{H^2} +\|\nabla_x \mathsf{u}\|_{H^2}\| \mathsf{B} \|_{H^2}  \big)\|f\|_{H_{x,v}^2}^2\\
&+C \|f\|_{H_{x,v}^2}^2+C\|\nabla_x \mathsf{B} \|_{H^2}^2\|f\|_{H_{x,v}^2}^2\\
\leq\,& C\big(1+\|\nabla_x \mathsf{u}\|_{H^2}^2+\|\nabla_x \mathsf{B} \|_{H^2}^2\big)\|f\|_{H_{x,v}^2}^2,
\end{align*}
where we have utilized the H\"{o}lder’s inequality and Lemma \ref{L1}.

Finally, we estimate the space-velocity mixed 
derivatives of $f$.  Select $\alpha$ and $\beta$ with $|\alpha|+|\beta|\leq 2$.
Applying $\partial_{\beta}^\alpha$ to \eqref{I3}, one has
\begin{align*}
&\partial_{\beta}^\alpha \partial_tf+\sum_{\beta^\prime\neq 0}C_{\beta}^{\beta^\prime}
\partial_{\beta^\prime}v\cdot\nabla_v\partial_{\beta-\beta^\prime}^{\alpha}f-\sum_{\alpha^\prime\neq 0}C_{\alpha}^{\alpha^\prime}\partial^{\alpha^\prime}(\mathsf{u}\times \mathsf{B} )\cdot\nabla_v\partial_{\beta}^{\alpha-\alpha^\prime}f\\
&\qquad =(\mathsf{u}\times \mathsf{B} )\cdot \nabla_v \partial_{\beta}^\alpha f-\sum_{\beta^\prime\neq 0}C_{\alpha}^{\alpha^\prime}C_{\beta}^{\beta^\prime}\partial_{\beta^\prime}v\times\partial^{\alpha^\prime} \mathsf{B} \cdot\nabla_v\partial_{\beta-\beta^{\prime}}^{\alpha-\alpha^\prime}f .  
\end{align*}
Multiplying the above equation by  $\partial_{\beta}^\alpha f$  
and integrating the result over $\mathbb{R}^3\times\mathbb{R}^3$, one gets
\begin{align*}
\frac{1}{2}\frac{\mathrm{d}}{\mathrm{d}t}\|\partial_{\beta}^\alpha f\|^2
=\,&-\sum_{\beta^\prime\neq 0}C_{\beta}^{\beta^\prime}\int_{\mathbb{R}^3\times\mathbb{R}^3}
(\partial_{\beta^\prime}v\cdot\nabla_v\partial_{\beta-\beta^\prime}^{\alpha}f)\partial_{\beta}^\alpha f\mathrm{d}x\mathrm{d}v\\
&+\sum_{\alpha^\prime\neq 0}C_{\alpha}^{\alpha^\prime}\int_{\mathbb{R}^3\times\mathbb{R}^3}\big(\partial^{\alpha^\prime}(\mathsf{u}\times \mathsf{B} )\cdot\nabla_v\partial_{\beta}^{\alpha-\alpha^\prime}f\big)\partial_{\beta}^\alpha f\mathrm{d}x\mathrm{d}v\\
&-\sum_{\beta^\prime\neq 0}C_{\alpha}^{\alpha^\prime}
C_{\beta}^{\beta^\prime}\int_{\mathbb{R}^3\times\mathbb{R}^3}
\big(\partial_{\beta^\prime}v\times\partial^{\alpha^\prime} \mathsf{B} \cdot\nabla_v\partial_{\beta-\beta^{\prime}}^{\alpha-\alpha^\prime}f\big)\partial_{\beta}^\alpha f\mathrm{d}x\mathrm{d}v\\
\leq\,&C\|f\|_{H_{x,v}^2}^2+C\big(\|\nabla_x \mathsf{B} \|_{H^2}\| \mathsf{u}\|_{H^2} +\|\nabla_x \mathsf{u}\|_{H^2}\| \mathsf{B} \|_{H^2}  \big)\|f\|_{H_{x,v}^2}^2\\
&+C \|f\|_{H_{x,v}^2}^2+C\|\nabla_x \mathsf{B} \|_{H^2}^2\|f\|_{H_{x,v}^2}^2\\
\leq\,& C\big(1+\|\nabla_x \mathsf{u}\|_{H^2}^2+\|\nabla_x \mathsf{B} \|_{H^2}^2\big)\|f\|_{H_{x,v}^2}^2.
\end{align*}
Combining the above estimates and taking summation,
we then obtain the desired estimate \eqref{G2.7}. 
\end{proof}
Using Gronwall's inequality for \eqref{G2.7} and the smallness of $\|f_0\|_{H_{x,v}^2}$, 
we obtain that
\begin{align}
\|f\|_{H_{x,v}^2}^2
\leq\,& e^{CT+C\int_0^t \big(\|\nabla_x \mathsf{u}(s)\|_{H^2}^2+\|\nabla_x \mathsf{B} (s)\|_{H^2}^2\big)\mathrm{d} s}\|f_0\|_{H_{x,v}^2}^2\nonumber\\
\leq\,&Ce^{CT}\|f_0\|_{H_{x,v}^2}^2\nonumber\\
\leq\,&C \|f_0\|_{H_{x,v}^2} ,
\end{align}
which implies
\begin{align}\label{G2.20}
\|f\|_{H_{x,v}^2}\leq\,&C \|f_0\|_{H_{x,v}^2}^{\frac{1}{2}}<\epsilon_1^{\frac{1}{3}}.
\end{align}
Thus we obtain the first inequality in \eqref{t2333}.

By a classical DiPerna-Lions theory \cite{DL-IM-1989} for the transport
equation with  our supposition on the regularity of $\mathsf{u}$ and $\mathsf{B} $, we finally achieve
the existence and uniqueness of $f$.
We notice that the existence and uniqueness of linear parabolic equations
\eqref{H2.1}$_2$--\eqref{H2.1}$_3$ can be easily achieved through the
semigroup theory (see\cite{kato-1973}). 
 Below we derive some estimates of   $u$ and $B$ by utilizing the  properties of $f$.


\subsubsection{The estimates of $u$ and $B$ to the system \eqref{H2.1}--\eqref{I--33} }

\begin{lem}\label{L2.3}
For the chosen $T$ stated in Lemma \ref{N2.3}, there exists a positive constant $\lambda_1>0$, such that
\begin{align}
&\frac{1}{2} \frac{\mathrm{d}}{\mathrm{d} t}\left(\|u\|^2+\|B\|^2\right)+\lambda_1\left(\|\nabla_x u\|^2+\|\nabla_x B\|^2\right)\nonumber \\ \label{H2.7}
&\qquad \quad \leq C\big((1+T)\epsilon_1^{\frac{1}{3}}+\epsilon_2\big)\left(\|\nabla_x u\|^2+\|\nabla_x B\|^2\right)
\end{align}
holds for all $0\leq t\leq T$.
\end{lem}
\begin{proof}
Multiplying the equations \eqref{H2.1}$_2$--\eqref{H2.1}$_3$ by $u$ and $B$, respectively, 
and   taking the results over  $\mathbb{R}^3 $ and summation, we obtain
\begin{align}
&\frac{1}{2} \frac{\mathrm{d}}{\mathrm{d} t}\left(\|u\|^2+\|B\|^2\right)+\|\nabla_x (u,B)\|^2\nonumber\\
=\,&-\int_{\mathbb{R}^3} \mathsf{u}\cdot\nabla_x u\cdot u\mathrm{d}x-\int_{\mathbb{R}^3} \mathsf{u}\cdot\nabla_x B\cdot B\mathrm{d}x+\int_{\mathbb{R}^3} \mathsf{B} \cdot\nabla_x u\cdot B\mathrm{d}x\nonumber\\
&+\int_{\mathbb{R}^3} \mathsf{B} \cdot\nabla_x B\cdot u\mathrm{d}x-\frac{1}{2}\sum_{i,j=1}^3\int_{\mathbb{R}^3} \mathsf{B} _i\partial_j B_i u_j\mathrm{d}x\nonumber\\
&+\int_{\mathbb{R}^3\times\mathbb{R}^3} 
u\cdot(u\times \mathsf{B} )f\mathrm{d}v\mathrm{d}x-\int_{\mathbb{R}^3\times\mathbb{R}^3} 
u\cdot(v\times B)f\mathrm{d}v\mathrm{d}x\nonumber\\
:=\,&\sum_{j=1}^{7}\mathcal{I}_{j}.
\end{align}
For the terms $\mathcal{I}_j$ ($j=1,\dots,6$), 
using the inequalities in Lemma \ref{L1}, we have
\begin{align*}&\mathcal{I}_1\leq C\|\mathsf{u}\|_{L^3}\|\nabla_x u\|\|u\|_{L^6}\leq C\|\mathsf{u}\|_{H^1}\|\nabla_x u\|^2\leq C \epsilon_1\|\nabla_x u\|^2,\\
&\mathcal{I}_2\leq C\|\mathsf{u}\|_{L^3}\|\nabla_x B\|\|B\|_{L^6}\leq C\|\mathsf{u}\|_{H^1}\|\nabla_x B\|^2\leq C \epsilon_1\|\nabla_x B\|^2,\\
&\mathcal{I}_3\leq C\|\mathsf{B} \|_{L^3}\|\nabla_x u\|\|B\|_{L^6}\leq C\|\mathsf{B} \|_{H^1}\|\nabla_x (u,B)\|^2\leq C \epsilon_1\|\nabla_x (u,B)\|^2,\\
&\mathcal{I}_4\leq C\|\mathsf{u}\|_{L^3}\|\nabla_x B\|\|B\|_{L^6}\leq C\|\mathsf{u}\|_{H^1}\|\nabla_x B\|^2\leq C \epsilon_1\|\nabla_x B\|^2,\\
&\mathcal{I}_5\leq C\|\mathsf{u}\|_{L^3}\|\nabla_x B\|\|B\|_{L^6}\leq C\|\mathsf{u}\|_{H^1}\|\nabla_x B\|^2\leq C \epsilon_1\|\nabla_x B\|^2,\\
&\mathcal{I}_6\leq C\|\mathsf{B} \|_{L^6}\|f\|_{H_{x,v}^2}\|u\|_{L^6}\|u\|_{L^6}\leq C\|\mathsf{B} \|_{H^1}\|\nabla_x u\|^2\leq C \epsilon_1\|\nabla_x u\|^2 .
\end{align*}
For the term $\mathcal{I}_7$, using H\"{o}lder’s inequality and \eqref{H2.6}, we have
\begin{align*}
\mathcal{I}_7\,&\leq C\int_{\mathcal{B}_x\times \mathcal{B}_v}|u||B||f|_2\mathrm{d}v\mathrm{d}x\ \\
&\leq C\|u\|_{L^6}\|B\|_{L^6}\|1\cdot {\chi({\mathcal{B}_x})}\|_{L^6(\mathcal{B}_x)}\|f\|_{H_{x,v}^2}    \\
&\leq C\|f\|_{H_{x,v}^2}(1+T)\left(\|\nabla_x u\|^2+\|\nabla_x B\|^2 \right)\\
&\leq C\epsilon_1^{\frac{1}{3}} (1+T)\big(\|\nabla_x u\|^2+\|\nabla_x B\|^2 \big),
\end{align*}
where $\chi(\cdot)$ is a characteristic   function. Collecting the above estimates together, we easily get the estimate \eqref{H2.7} 
and hence  complete the proof of Lemma \ref{L2.3}.
\end{proof}
\begin{lem}\label{L2.4}
For the chosen $T$ stated in Lemma \ref{N2.3}, there exists a positive constant $\lambda_2>0$, such that
\begin{align}
&\frac{1}{2} \frac{\mathrm{d}}{\mathrm{d} t} \sum_{1 \leq|\alpha| \leq 2}\left(\|\partial^\alpha u\|^2+\|\partial^\alpha B\|^2\right) +\lambda_2 \sum_{1 \leq|\alpha| \leq 2}
\Big(\|\nabla_x \partial^\alpha u\|^2 +\|\nabla_x \partial^\alpha B\|^2 \Big) \nonumber\\
&\qquad \qquad
\leq C(\epsilon_1^{\frac{1}{3}}+\epsilon_2)\left(\|\nabla_x u\|_{H^2}^2+\|\nabla_x B\|_{H^2}^2\right) \label{H2.9}
\end{align}
holds for all $0\leq t\leq T$.
\end{lem}
\begin{proof}
Applying $\partial^{\alpha}$ with $1 \leq|\alpha| \leq 2$ 
to the equations \eqref{H2.1}$_2$--\eqref{H2.1}$_3$,  
and multiplying the results by $\partial^\alpha u$ and $ \partial^\alpha B$, respectively, then taking the integration 
over $\mathbb{R}^3$ and summing   them together, we   obtain
\begin{align}
& \frac{1}{2} \frac{\mathrm{d}}{\mathrm{d} t} \left(\left\|\partial^\alpha u\right\|^2+\left\|\partial^\alpha B\right\|^2\right)+\|\nabla_x \partial^\alpha u\|^2+\|\nabla_x \partial^\alpha B\|^2 \nonumber\\
=\,&-\int_{\mathbb{R}^3}\partial^{\alpha}(\mathsf{u}\cdot \nabla_x B)\cdot\partial^{\alpha}B\mathrm{d}x+\int_{\mathbb{R}^3}\partial^{\alpha}(\mathsf{B} \cdot \nabla_x B)\cdot\partial^{\alpha}u\mathrm{d}x+\int_{\mathbb{R}^3}\partial^{\alpha}(\mathsf{B} \cdot \nabla_x u)\cdot\partial^{\alpha}B\mathrm{d}x\nonumber\\
&-\frac{1}{2}\sum_{i,j=1}^3\int_{\mathbb{R}^3}\partial^{\alpha}(\mathsf{B} _i\partial_j B_i)\partial^{\alpha}u_j\mathrm{d}x-\int_{\mathbb{R}^3}\partial^{\alpha}(\mathsf{u}\cdot \nabla_x u)\cdot\partial^{\alpha}u\mathrm{d}x\nonumber\\
&+\int_{\mathbb{R}^3\times\mathbb{R}^3} 
\partial^{\alpha}u\cdot\partial^{\alpha}\big((u\times \mathsf{B} )f\big)\mathrm{d}v\mathrm{d}x-\int_{\mathbb{R}^3\times\mathbb{R}^3} 
\partial^{\alpha}u\cdot\partial^{\alpha}\big((v\times B)f\big)\mathrm{d}v\mathrm{d}x\nonumber\\
:=\,&\sum_{j=8}^{14}\mathcal{I}_{j} .
\end{align}
With the aid of Lemma \ref{L1}, we have
\begin{align*}
\mathcal{I}_{8}\leq C\|\partial^{\alpha}(\mathsf{u}\cdot\nabla_x B)\|\|\partial^{\alpha} B\| 
 \leq C\|\mathsf{u}\|_{H^2}\|\nabla_x B\|_{H^2}^2 
 \leq C\epsilon_2\|\nabla_x B\|_{H^2}^2.
\end{align*}
Similarly, we   obtain
\begin{align*}
\mathcal{I}_{9}&\leq C\|\partial^{\alpha}(\mathsf{B} \cdot\nabla_x B)\|\|\partial^{\alpha} u\|\\
&\leq C\|\mathsf{B} \|_{H^2}\left(\|\nabla_x B\|_{H^2}^2+\|\nabla_x u\|_{H^2}^2\right) \\
&\leq C\epsilon_2\left(\|\nabla_x B\|_{H^2}^2+\|\nabla_x u\|_{H^2}^2\right),\\
\mathcal{I}_{10}&\leq C\|\partial^{\alpha}(\mathsf{B} \cdot\nabla_x u)\|\|\partial^{\alpha} B\|\\
&\leq C\|\mathsf{B} \|_{H^2}\left(\|\nabla_x B\|_{H^2}^2+\|\nabla_x u\|_{H^2}^2\right)\\
&\leq C\epsilon_2\left(\|\nabla_x B\|_{H^2}^2+\|\nabla_x u\|_{H^2}^2\right),\\
\mathcal{I}_{11}&\leq C\|\partial^{\alpha}(\mathsf{B} _i\partial_j B_i)\|\|\partial^{\alpha} u_j\|\\
&\leq C\|\mathsf{B} \|_{H^2}\left(\|\nabla_x B\|_{H^2}^2+\|\nabla_x u\|_{H^2}^2\right)\\
&\leq C\epsilon_2\left(\|\nabla_x B\|_{H^2}^2+\|\nabla_x u\|_{H^2}^2\right),\\
\mathcal{I}_{12}&\leq C\|\partial^{\alpha}(\mathsf{u}\cdot\nabla_x u)\|\|\partial^{\alpha} u\|\\
&\leq C\|\mathsf{u}\|_{H^2}\|\nabla_x u\|_{H^2}^2\\
&\leq C\epsilon_2\|\nabla_x u\|_{H^2}^2,\\
\mathcal{I}_{13}&\leq C\|\mathsf{B} \|_{H^2}\|f\|_{H_{x,v}^2}\|\nabla_x u\|_{H^2}^2\\
&\leq C\|\mathsf{B} \|_{H^2}\|\nabla_x u\|_{H^2}^2\\
&\leq C\epsilon_2\|\nabla_x u\|_{H^2}^2.
\end{align*}
For the remaining term $\mathcal{I}_{14}$, 
using H\"{o}lder’s inequality and Lemma \ref{L1}, 
we have
\begin{align*}
\mathcal{I}_{14}\,&\leq C\int_{\mathbb{R}^3\times\mathbb{R}^3} 
|\partial^{\alpha}u||\partial^{\alpha}\big(( v\times B)f\big)|\mathrm{d}v\mathrm{d}x\\
&\leq C\|f\|_{H_{x,v}^2}\left(\|B\|_{L^{\infty}}\|\nabla_x u\|_{H^2}+\|u\|_{L^{\infty}}\|\nabla_x B\|_{H^2}\right)
+C\|f\|_{H_{x,v}^2}\left(\|\nabla_x u\|_{H^2}^2+\|\nabla_x B\|_{H^2}^2\right)\\
&\leq C\|f\|_{H_{x,v}^2}\|\nabla_x B\|_{H^2}\|\nabla_x u\|_{H^2}
+C\|f\|_{H_{x,v}^2}\left(\|\nabla_x u\|_{H^2}^2+\|\nabla_x B\|_{H^2}^2\right)\\
&\leq C\epsilon_1^{\frac{1}{3}}\left(\|\nabla_x u\|_{H^2}^2+\|\nabla_x B\|_{H^2}^2\right).
\end{align*}
Combining the above estimates together, and then taking summation over $1 \leq|\alpha| \leq 2$, we obtain the estimate \eqref{H2.9}. Thus, we complete the proof of Lemma \ref{L2.4}.
\end{proof}

Define the temporal energy functional 
$\mathcal{E}(t)$ and the corresponding dissipation rate
functional $\mathcal{D}(t)$ as follows:
\begin{align}\nonumber
\mathcal{E}(t) =\|u\|_{H^2}^2+\|B\|_{H^2}^2,\ \
\mathcal{D}(t) =\|\nabla_x u\|_{H^2}^2+\|\nabla_x B\|_{H^2}^2.
\end{align}
Adding \eqref{H2.7} to \eqref{H2.9}, we have 
\begin{align}\nonumber %
\frac{\mathrm{d}}{\mathrm{d} t}\mathcal{E}(t)+\lambda_3\mathcal{D}(t)
\leq C\big((1+T)\epsilon_1^{\frac{1}{3}}+\epsilon_2\big)\mathcal{D}(t) 
\end{align}
for some $\lambda_3>0$.
Using the smallness of $\epsilon_1$ and $\epsilon_2$, we deduce that
\begin{align}\label{H2.14}
\frac{\mathrm{d}}{\mathrm{d} t}\mathcal{E}(t)+\lambda_4\mathcal{D}(t)\leq 0 
\end{align}
for some $\lambda_4>0$.
Integrating \eqref{H2.14} over $(0,t)$ yields
\begin{align}\nonumber 
\mathcal{E}(t)+\lambda_4\int_0^t \mathcal{D}(s)\mathrm{d} s\leq \mathcal{E}(0)   
\end{align}
for all $0\leq t \leq T$.

Thus, we have
\begin{gather*}
\|(u,B)\|_{H^2}\leq \|(u_0,B_0)\|_{H^2}<  \epsilon_2,\\
\int_0^t\|\nabla_x u(s)\|_{H^2}^2+\|\nabla_x B(s)\|_{H^2}^2\mathrm{d}s\leq  \mathcal{E}(0) 
\end{gather*}
for all $0\leq t \leq T$,  which implies that 
\begin{align*}
\|(u,B)\|_{C ([0,T];H^{2})}<\epsilon_2,\quad 
\|(u,B)\|_{L^2 (0,T;H^{3})}<\epsilon_2.    
\end{align*}
Therefore, we obtain the second and third inequalities in \eqref{t2333} and complete the proof of Theorem \ref{reaaa}.

\subsection{Global existence of strong solutions to the problem  \eqref{I3}--\eqref{I--3}} 
In this subsection, we will show the global existence and uniqueness of strong solutions to the problem \eqref{I3}--\eqref{I--3}. 
To this end, we first  construct the approximation 
sequences $\{(f^{n+1},u^{n+1},B^{n+1})\}$ satisfying
\begin{equation}\label{G2.22}\left\{
\begin{aligned}
&\partial_t f^{n+1}+v\cdot\nabla_x f^{n+1}+\big((v-u^{n})\times B^{n}\big)\cdot\nabla_v f^{n+1} =0, \\
&\partial_t u^{n+1}+u^{n}\cdot\nabla_x u^{n+1}+\nabla_x P^{n+1}-B^{n}\cdot\nabla_x B^{n+1}-\Delta_x u^{n+1}+\frac{1}{2}\sum_{i,j=1}^3 B_i^n \partial_j B_i^{n+1} \\
&\qquad  =\int_{\mathbb{R}^3}(u^{n+1}\times B^{n})f^{n+1}\mathrm{d}v-\int_{\mathbb{R}^3}(v\times B^{n+1})f^{n+1}\mathrm{d}v  \\
& \partial_t B^{n+1}-\Delta_x B^{n+1}+u^n\cdot\nabla_x B^{n+1}-B^n\cdot\nabla_x u^{n+1}=0,\\
&\nabla\cdot u^{n+1}=0,\quad \nabla\cdot B^{n+1}=0,
\end{aligned}\right.
\end{equation}
 with the initial data:
\begin{align}\label{G2.23}
(f^{n+1}(t,x,v),u^{n+1}(t,x),B^{n+1}(t,x))|_{t=0}=(f_0(x,v),u_0(x),B_0(x)),
\end{align}
for all $n\geq 0$ and $(x,v)\in\mathbb{R}^3\times\mathbb{R}^3$. In addition,
\begin{align}\label{G2.24}
(f^0(t,x,v),u^0(t,x),B^0(t,x))\equiv(f_0(x,v),u_0(x),B_0(x)),\quad(t,x,v)\in [0,T]\times\mathbb{R}^3\times\mathbb{R}^3.
\end{align}

As a direct application of Theorem \ref{reaaa}, 
the existence and uniqueness of solutions to 
the problem \eqref{G2.22}--\eqref{G2.24} 
can be guaranteed.
\begin{prop}\label{L2.7}
Suppose that $f_0(x,v)$ has a compact
support in $x$ and $v$,
assume that there
exists $\epsilon_1=\epsilon_1(T)>0$ and $\epsilon_2=\epsilon_2(T)>0$, such that
if $\|f_0\|_{H_{x,v}^2}<\epsilon_1$ and
$\|(u_0,B_0)\|_{H^{2}}<\epsilon_2$, 
then for each $n\geq 0$, 
there exists a unique solution $(f^{n+1}, u^{n+1}, B^{n+1})$ to the problem \eqref{G2.22}--\eqref{G2.24} such that
\begin{align}\label{n+1}
\|f^{n+1}\|_{C([0,T];H_{x,v}^2)}< \epsilon_1^{\frac{1}{3}},\quad\|(u^{n+1},B^{n+1})\|_{C ([0,T];H^{2})}<\epsilon_2,\quad 
\|(u^{n+1},B^{n+1})\|_{L^2 (0,T;H^{3})}<\epsilon_2.
\end{align}
\end{prop}

\begin{rem} 
In Proposition \ref{L2.7}, we can choose  $0<\epsilon_1\ll 1$  with an order of $e^{-\mathcal{O}(1)T}$,
and $0<\epsilon_2\ll 1$ with an order of $\frac{\mathcal{O}(1)}{T+1}$. 
\end{rem}
We now show that the approximation sequence $\{(f^n,u^n,B^n)\}_{n\geq 1}$ is a Cauchy sequence in 
$C([0, T]; L^2(\mathbb{R}^3\times\mathbb{R}^3))\times C([0, T]; H^1(\mathbb{R}^3))\cap L^2(0,T;H^2(\mathbb{R}^3)) \times C([0, T]; H^1(\mathbb{R}^3)) \cap L^2(0,T; \linebreak H^2(\mathbb{R}^3))$.
By a direct calculation,
$(f^{n+1}-f^n,u^{n+1}-u^n,B^{n+1}-B^n)$ satisfies the 
following equations:

\begin{equation}\label{G2.33}
\left\{
\begin{aligned}
&\partial_t (f^{n+1}-f^n)+v\cdot\nabla_x (f^{n+1}-f^n)+v\times( B^{n}-B^{n-1})\cdot\nabla_v f^{n}  \\
&\qquad +v\times B^n\cdot \nabla_v(f^{n+1}-f^n)-(u^n\times B^n)\cdot\nabla_v (f^{n+1}-f^n)\\
&\qquad -(u^{n}-u^{n-1})\times B^n\cdot \nabla_v f^{n}-u^{n-1}\times(B^{n}-B^{n-1})\cdot\nabla_v f^{n}=0,      \\  
&\partial_t  (u^{n+1}-u^n)+u^{n}\cdot\nabla_x (u^{n+1}-u^n)+(u^{n}-u^{n-1})\nabla_x u^n\\
&\qquad+\nabla_x(P^{n+1}-P^n)-(B^n-B^{n-1})\nabla_x B^{n+1}-B^{n-1}\cdot\nabla_x (B^{n+1}-B^n)\\
&\qquad-\Delta_x(u^{n+1}-u^n)+\frac{1}{2}\sum_{i,j=1}^3 B_i^{n-1} \partial_j (B_i^{n+1}-B_i^{n})+\frac{1}{2}\sum_{i,j=1}^3 (B_i^n-B_i^{n-1}) \partial_j B_i^{n+1}                    \\
&\quad=-\int_{\mathbb{R}^3}v\times(B^{n+1}-B^{n})f^{n}\mathrm{d}v-\int_{\mathbb{R}^3}v\times B^{n+1}(f^{n+1}-f^{n})\mathrm{d}v\\
&\qquad+\int_{\mathbb{R}^3}(u^{n+1}\times B^n)(f^{n+1}-f^{n})\mathrm{d}v
+\int_{\mathbb{R}^3}(u^{n+1}-u^n)\times B^n f^{n}\mathrm{d}v\\
&\qquad-\int_{\mathbb{R}^3}u^n\times(B^{n}-B^{n-1})f^{n}\mathrm{d}v,                        \\
& \partial_t (B^{n+1}-B^n)-\Delta_x (B^{n+1}-B^n)+u^n\cdot\nabla_x (B^{n+1}-B^n)+\nabla_x B^n\cdot(u^n-u^{n-1})\\ 
&\qquad -(B^n-B^{n-1})\nabla_x u^{n+1}- B^{n-1}\cdot\nabla_x(u^{n+1}-u^{n})=0,\\  
&\nabla\cdot (u^{n+1}-u^n)=0,\quad \nabla\cdot (B^{n+1}-B^n)=0.
\end{aligned}\right.
\end{equation}
Multiplying \eqref{G2.33}$_1$ by $f^{n+1}-f^n$ and integrating the results over $\mathbb{R}^3\times\mathbb{R}^3$, one has
\begin{align*}
\frac{1}{2}\frac{\mathrm{d}}{\mathrm{d}t}\|f^{n+1}-f^n\|_{L^2}^2
=&-\int_{\mathbb{R}^3}v\times(B^{n+1}-B^n)\cdot\nabla_v f^{n}(f^{n+1}-f^n) \mathrm{d}v\mathrm{d}x\\
&+\int_{\mathbb{R}^3}(u^{n}-u^{n-1})\times B^n\cdot\nabla_v f^{n}(f^{n+1}-f^n) \mathrm{d}v\mathrm{d}x\\
&+\int_{\mathbb{R}^3}u^{n-1}\times (B^n-B^{n-1})\cdot\nabla_v f^{n}(f^{n+1}-f^n) \mathrm{d}v\mathrm{d}x\\
\leq \,& C\|f^{n+1}\|_{H_{x,v}^2}\|B^{n}-B^{n-1}\|_{L^6}\|f^{n+1}-f^n\|_{L^2}\\
&+C\|f^{n}\|_{H_{x,v}^2}\|B^n\|_{L^6}\|u^{n+1}-u^n\|_{L^6}\|f^{n+1}-f^n\|_{L^2}\\
&+C\|f^{n}\|_{H_{x,v}^2}\|u^{n-1}\|_{L^6}\|B^{n}-B^{n-1}\|_{L^6}\|f^{n+1}-f^n\|_{L^2}\\
\leq \,& C\epsilon_1^{\frac{1}{3}}\big(\|\nabla_x(B^{n}-B^{n-1})\|_{L^2}^2+\|\nabla_x(u^{n}-u^{n-1})\|_{L^2}^2\big.\\
\big.&+\|f^{n+1}-f^n\|_{L^2}^2  \big).
\end{align*}
Using the same method as the above arguments, we  obtain
\begin{align*}
&\frac{1}{2}\frac{\mathrm{d}}{\mathrm{d}t}\big\{\|u^{n+1}-u^n\|_{H^1}^2+\|B^{n+1}-B^n\|_{H^1}^2\big\}\\
&+\lambda_{5}\big(\|\nabla_x (u^{n+1}-u^n)\|_{H^1}^2+\|\nabla_x(B^{n+1}-B^n)\|_{H^1}^2\big)\\
\leq\,& C(\epsilon_1^{\frac{1}{3}}+\epsilon_2)
\big(\|u^{n+1}-u^n\|_{H^1}^2+\|B^{n+1}-B^n\|_{H^1}^2\big)+C\epsilon_2\|f^{n+1}-f^{n}\|_{L^2}^2\\
&+C(\epsilon_1^{\frac{1}{3}}+\epsilon_2)\big(\|u^{n}-u^{n-1}\|_{H^1}^2+\|B^{n}-B^{n-1}\|_{H^1}^2\big) 
\end{align*}
for some $\lambda_{5}>0$.

Therefore, we have
\begin{align}
&\frac{1}{2}\frac{\mathrm{d}}{\mathrm{d}t}\big\{\|f^{n+1}-f^{n}\|_{L^2}^2+\|u^{n+1}-u^n\|_{H^1}^2
+\|B^{n+1}-B^n\|_{H^1}^2\big\}\nonumber\\
&+\lambda_{5}\big(\|\nabla_x (u^{n+1}-u^n)\|_{H^1}^2+\|\nabla_x(B^{n+1}-B^n)\|_{H^1}^2\big)\nonumber\\
\leq\,& C(\epsilon_1^{\frac{1}{3}}+\epsilon_2)\big(\|f^{n+1}-f^{n}\|_{L^2}^2+\|u^{n+1}-u^n\|_{H^1}^2
+\|B^{n+1}-B^n\|_{H^1}^2\big)\nonumber\\\label{GG2.28}
&+C(\epsilon_1^{\frac{1}{3}}+\epsilon_2)\big(\|f^{n}-f^{n-1}\|_{L^2}^2+\|u^{n}-u^{n-1}\|_{H^1}^2
+\|B^{n}-B^{n-1}\|_{H^1}^2\big).
\end{align}

Applying  the Gronwall's lemma (see Lemma 3 in \cite{BDGM-die-2009}) to the above inequality
 gives  

\begin{align*}
&\|f^{n+1}-f^{n}\|_{C([0,T];L^2)}+\|u^{n+1}-u^{n}\|_{C ([0,T];H^{1})}+\|B^{n+1}-B^{n}\|_{C ([0,T];H^{1})}< \frac{C(T)^{n+1}}{n!},
\end{align*}
which implies that 
\begin{align*}
&\|\nabla_x (u^{n+1}-u^{n})\|_{L^2 (0,T;H^{1})}+\|\nabla_x (B^{n+1}-B^{n})\|_{L^2 (0,T;H^{1})}< \frac{C(T)^{n+1}}{n!}    
\end{align*}
for all $n\geq 0$.

Therefore, we  conclude that $\{(f^n,u^n,B^n)\}_{n\geq 0}$ is a Cauchy sequence in the Banach space
$C([0, T]; L^2(\mathbb{R}^3\times\mathbb{R}^3))\times C([0, T]; H^1(\mathbb{R}^3))\cap L^2(0,T;H^2(\mathbb{R}^3)) \times C([0, T]; H^1(\mathbb{R}^3)) \cap L^2(0,T; \linebreak H^2(\mathbb{R}^3))$.
And there exist the 
limit functions $(f,u,B)$ such that
\begin{align}
f^{n}\to f\quad\mathrm{in}&\quad C([0,T];L^{2}(\mathbb{R}^3\times\mathbb{R}^{3})),\label{G2.38-1}\\
\quad B^{n}\to B\quad\mathrm{in}&\quad C([0,T];H^{1}(\mathbb{R}^3))\cap L^2(0,T;H^2(\mathbb{R}^3)),\label{G2.38-2}\\\label{G2.38}
u^{n}\to u\quad\mathrm{in}&\quad C([0,T];H^{1}(\mathbb{R}^3))\cap L^2(0,T;H^2(\mathbb{R}^3)).
\end{align}
By the Gagliardo-Nirenberg interpolation inequality
\begin{align*}
\|f\|_{H^1_{x,v}}\leq \|f\|_{L^2}^{\frac{1}{2}}\|f\|_{H^2_{x,v}}^{\frac{1}{2}},    
\end{align*}
we know that 
\begin{align}\label{G2.39}
f^n\to f\quad\mathrm{in}\quad C([0,T];H^1_{x,v}(\mathbb{R}^3\times\mathbb{R}^3)).
\end{align}

We now claim that $(f,u,B)\in C([0,T];H^2_{x,v}(\mathbb{R}^3\times\mathbb{R}^3))\times C([0,T];H^2(\mathbb{R}^3))\times\ C([0,T];\linebreak H^2(\mathbb{R}^3))$.  In fact, it follows from Proposition \ref{L2.7} 
that there exists a convergent subsequence $(f^{n_k},u^{n_k},  B^{n_k})$ such that $f^{n_k}\stackrel{*}{\rightharpoonup}\mathfrak{f}$, 
$ (u^{n_{k}}, B^{n_{k}})\rightharpoonup(\mathfrak{u}, \mathfrak{B})$, and
 \begin{align}
\|\mathfrak{f}(t)\|_{H_{x,v}^{2}}\leq \,&\operatorname*{lim}_{k\to\infty}\inf\|f^{n_{k}}(t)\|_{H_{x,v}^{2}},\label{G2.40-1}\\
 \|\mathfrak{u}(t)\|_{H^{2}}\leq \,& \operatorname*{lim}_{k\to\infty}\inf\|u^{n_{k}}(t)\|_{H^{2}},\label{G2.40-2}\\
 \|\mathfrak{B}(t)\|_{H^{2}}\leq \,& \operatorname*{lim}_{k\to\infty}\inf\|B^{n_{k}}(t)\|_{H^{2}} \label{G2.40}
\end{align}
for any $0\leq t\leq T$. 

For $(\mathfrak{f},\mathfrak{u},\mathfrak{B})\in L^\infty([0,T];H^2_{x,v}(\mathbb{R}^3\times\mathbb{R}^3))\times L^\infty([0,T];H^2(\mathbb{R}^3))\times L^\infty([0,T];H^2(\mathbb{R}^3))$,
combining the above fact with 
\eqref{G2.38}--\eqref{G2.40}, we   infer that
$(\mathfrak{f},\mathfrak{u},\mathfrak{B}) = (f,u,B) $ in $L^\infty([0,T];\linebreak H^2_{x,v}(\mathbb{R}^3\times\mathbb{R}^3))\times L^\infty([0,T]; H^2(\mathbb{R}^3))\times L^\infty([0,T];H^2(\mathbb{R}^3)) $.

Additionally, we  can show that $(f,u,B)\in C_w([0,T];H^2_{x,v}(\mathbb{R}^3\times\mathbb{R}^3))\times C_w([0,T];H^2(\mathbb{R}^3))\times C_w([0,T];H^2(\mathbb{R}^3))$ by taking arguments similar to the above.
For brevity, below we only show 
that $B\in C_w([0,T]; $   $H^2(\mathbb{R}^3))$,
and  we shall use the method  {introduced} in \cite{BV-2022} (See Chapter 2 in \cite{BV-2022}).
Notice that
\begin{align*}
\Big|\|B^{n+1}(t)\|_{L^2}^2-\|B^{n+1}(s)\|_{L^2}^2\Big|
=\,& \Big|\int_s^t\frac{\mathrm{d}}{\mathrm{d}\eta}\|B^{n+1}(\eta)\|_{L^2}^2\mathrm{d}\eta\Big| \\
\leq\,& C|t-s|\left(\|\nabla_x u^{n+1}\|^2+\|\nabla_x B^{n+1}\|^2\right)\\
\leq\,& C\delta_2^2|t-s|.
\end{align*}
The above fact  implies that $B^{n+1}\in \text{Lip}([0,T];L^2(\mathbb{R}^3))$.

We point out that $\langle B(t),\varphi\rangle_{L^{2}}$ 
is a continuous function of time, 
where $\varphi\in$ $H^{-2}(\mathbb{R}^3))$ is a test functions.  
Let's choose  $\varphi\in H^{-2}(\mathbb{R}^3))$ with the norm $\|\varphi\|_{H^{-2}}=1$. 
Let $\delta^{*}\in(0,1/4]$ be arbitrary. 
Thanks to the density of function space, we can find an element $\psi\in L^{2}$ 
such that
\begin{align*}
\|\varphi-\psi\|_{H^{-2}}\leq\delta^{*}.
\end{align*}

Due to \eqref{G2.38-2}, it holds that
\begin{align*}
\langle B^{n+1}(t),\psi\rangle_{{L^2}\times L^{2}}\to\langle B(t),\psi\rangle_{L^{2}\times L^{2}},
\end{align*}
as $n \to {\infty}$, uniformly in $t\in[0,T]$. 
Thus, we may find $n$ is sufficiently
large such that
\begin{align*}
\sup_{t\in[0,T]}|\langle B(t)-B^{n+1}(t),\psi\rangle|\leq\delta^{*}.
\end{align*}
Moreover, thanks to the
$B^{n+1}\in\text{Lip}([0,T];L^2(\mathbb{R}^3))$,
we know that
there exists $\tau^{*}>0$ such that
\begin{align*}
\sup\limits_{t,s\in[0,T];|t-s|\le\tau^{*}}\|B^{n+1}(t)-B^{n+1}(s)\|_{L^{2}}\le\delta^{*}(1+\|\psi\|_{L^{2}})^{-1}.
\end{align*}
From the above three ingredients, for $|t-s|\leq\tau^{*}$, we have
\begin{align*}
|\langle B(t)-B(s),\varphi\rangle|\leq\,&|\langle B^{n+1}(t)-B^{n+1}(s),\psi\rangle|+|\langle B(t)-B(s),\varphi-\psi\rangle|\\
&+|\langle B(t)-B^{n+1}(t),\psi\rangle|+|\langle B(s)-B^{n+1}(s),\psi\rangle|\\
\leq\,&\delta^{*}+2\delta^{*}\|B\|_{L^{\infty}([0,T];H^{2})}+2\delta^{*}.
\end{align*}
Since $\delta^{*}\in(0,1/4]$ is 
arbitrary, we   deduce that $B$ is weakly uniformly
continuous with values in $H^{2}$.

We now claim that for any $t_0\in [0,T)$, it holds that  {\begin{align*}
\|f( t) \|_{H_{x,v}^{2}}&\to \|f( t_{0}) \|_{H_{x,v}^{2}}, \\
\|u ( t) \|_{H^{2}}&\to \|u( t_{0}) \|_{H^{2}}, \\
\|B(t)\|_{H^2}&\to\|B(t_0)\|_{H^2}
\end{align*}
  as   $ t\to t_0^+$.}
For this, we   adopt the arguments used 
in  {Lemmas \ref{L2.0}--\ref{L2.4}} to obtain 
 {\begin{align}\label{G2.41}
\frac{\mathrm{d}}{\mathrm{d}t}\|f^n(t)\|_{H_{x,v}^2}^2\leq C, \quad\frac{\mathrm{d}}{\mathrm{d}t}\|B^n(t)\|_{H^2}^2\leq C,\quad \frac{\mathrm{d}}{\mathrm{d}t}\|u^n(t)\|_{H^2}^2\leq C,    
\end{align}}
where $(f^n,u^n,B^n)$ is the approximation solutions to the problem \eqref{G2.22}--\eqref{G2.24}. Without loss of
generality, we   suppose $t_0=0$. It follows from \eqref{G2.41} that 
{\begin{align*}
\|f^n(t)\|_{H_{x,v}^2}\leq\,& \|f_0\|_{H_{x,v}^2}+Ct,\\
 \|u^n(t)\|_{H^2}\leq\,& \|u_0\|_{H^2}+Ct,\\
 \|B^n(t)\|_{H^2}\leq \,&\|B_0\|_{H^2}+Ct.
\end{align*}}Then the lower semi-continuity and 
weak continuity explain   our claim is correct.
We  arrive at 
{\begin{align*}
\lim\limits_{t\to0+}\|f(t)-f_0\|_{H_{x,v}^2}=\lim\limits_{t\to0+}\|B(t)-B_0\|_{H^2}
=\lim\limits_{t\to0+}\|u(t)-u_0\|_{H^2}=0,
\end{align*}}
by applying  the facts that $f(t)\stackrel{*}{\rightharpoonup}f_0$ in $H^2_{x,v}(\mathbb{R}^3\times\mathbb{R}^3)$ and
$ (u^{n_{k}},B^{n_{k}})\rightharpoonup({u_0}, {B_0})$ in $H^2(\mathbb{R}^3)\times H^2(\mathbb{R}^3)$.
Considering the continuity from the left-hand side, we use a change of variable $t\mapsto T-t$ for $f$. 
However, the equations with regard to $u$ and $B$ 
are not time-reversible.
Therefore, we need to use 
the smoothing effect of diffusion  
for $u$ and $B$ to achieve the
left continuity (see Section 3.2 in \cite{MB-2002}).
Then, we complete the proof of existence of 
strong solutions.

In order to show the uniqueness of strong solutions, we set $({f_1},{u_1},{B_1})\in C([0,T];H^2_{x,v}(\mathbb{R}^3\times\mathbb{R}^3))\times C([0,T];H^2(\mathbb{R}^3))\times\ C([0,T];H^2(\mathbb{R}^3))$ be another solution to the Cauchy problem \eqref{I3}-\eqref{I--3}. Using similar argument to \eqref{GG2.28}, it follows that
\begin{align*}
&\|(f-{f_1})(t)\|_{L^{2}}^2+\|(u-{u_1})(t)\|_{H^{1}}^2+\|(B-{B_1})(t)\|_{H^{1}}^2\\
 &\qquad \leq \int_0^t\big(\|(f-{f_1})(s)\|^2_{L^{2}}+\|(u-{u_1})(s)\|^2_{H^{1}}+\|(B-{B_1})(s)\|^2_{H^{1}}\big)\mathrm{d}s,
\end{align*}
and
\begin{align*}
&\|(f-{f_1})(0)\|_{L^{2}}^2+\|(u-{u_1})(0)\|_{H^{1}}^2+\|(B-{B_1})(0)\|_{H^{1}}^2=0 
\end{align*}
for all $0\leq t \leq T$, which  implies that 
\begin{align*}
\|(f-{f_1})(t)\|_{L^{2}}^2+\|(u-{u_1})(t)\|_{H^{1}}^2+\|(B-{B_1})(t)\|_{H^{1}}^2=0    
\end{align*}
for all $t\in[0,T]$, and
\begin{align*}
f\equiv {f_1}&\quad \text{in} \quad C([0, T]; L^2(\mathbb{R}^3\times\mathbb{R}^3)),\\
u \equiv {u_1}&\quad \text{in} \quad C([0,T];H^{1}(\mathbb{R}^3))\cap L^2(0,T;H^2(\mathbb{R}^3)),\\
B \equiv {B_1}&\quad \text{in} \quad C([0,T];H^{1}(\mathbb{R}^3))\cap L^2(0,T;H^2(\mathbb{R}^3)).               
\end{align*}
Therefore, we obtain \eqref{t3333}. We further obtain \eqref{KKK} from \eqref{n+1} and 
\eqref{G2.40-1}--\eqref{G2.40}. Thus, we complete the existence part of Theorem \ref{I4}.  

\smallskip
\section{Decay rates of strong solutions to the problem \eqref{I3}--\eqref{I--3} in $\mathbb{R}^3$ } 
The preceding section gives the existence and uniqueness of strong solutions to the problem \eqref{I3}--\eqref{I--3}.    This section is devoted to investigating the time decay of 
strong solutions to it.


Let   $(f,u,B)$ be the strong solution to the Cauchy problem \eqref{I3}--\eqref{I--3} 
constructed in the above section on $0\leq t\leq T$, and  satisfy
\begin{align}
\sup_{ {0\leq t\leq T}}\|f\|_{H_{x,v}^2}\leq  \delta_1,\ \quad 
\sup_{ {0\leq t\leq T}}\|(u,B)\|_{H^{2}}\leq  \delta_2,  \label{HH2.1}
\end{align}
where $\delta_1$ has an order of $e^{-\mathcal{O}(1)T}$ and
$\delta_2$ has an order of $\frac{\mathcal{O}(1)}{T+1}$. Then, applying the arguments 
to Lemmas \ref{L2.0}--\ref{L2.4}, we know that $(f,u,B)$ satisfies 
\begin{align}
&\frac{1}{2} \frac{\mathrm{d}}{\mathrm{d} t}\|f\|_{H_{x,v}^2}^2
\leq C\big(1+\|\nabla_x u\|_{H^2}^2+\|\nabla_x B\|_{H^2}^2\big)\|f\|_{H_{x,v}^2}^2,\label{GG2.7}\\
&\frac{1}{2} \frac{\mathrm{d}}{\mathrm{d} t}\left(\|u\|^2+\|B\|^2\right)+\lambda_6\left(\|\nabla_x u\|^2+\|\nabla_x B\|^2\right)\nonumber \\
&\qquad \quad \leq C\big((1+T)\delta_1+\delta_2\big)\left(\|\nabla_x u\|^2+\|\nabla_x B\|^2\right),\label{HH2.7}\\
&\frac{1}{2} \frac{\mathrm{d}}{\mathrm{d} t} \sum_{1 \leq|\alpha| \leq 2}\left(\|\partial^\alpha u\|^2+\|\partial^\alpha B\|^2\right) +\lambda_7 \sum_{1 \leq|\alpha| \leq 2}
\Big(\|\nabla_x \partial^\alpha u\|^2 +\|\nabla_x \partial^\alpha B\|^2 \Big) \nonumber\\
&\qquad\leq C(\delta_1+\delta_2)\left(\|\nabla_x u\|_{H^2}^2+\|\nabla_x B\|_{H^2}^2\right) \label{HH2.9}
\end{align}
 for all $0\leq t\leq T$, and some $\lambda_6, \lambda_7>0$.
  
Define the temporal energy functional 
$\mathcal{E}_1(t)$ and the corresponding dissipation rate
functional $\mathcal{D}_1(t)$  to $(u,B)$ of the system \eqref{I3} as follows:
\begin{align}\label{GGG2.13}
&\mathcal{E}_1(t)=\|u\|^2+\|B\|^2+\sum_{1\leq |\alpha|\leq 2}\left(\|\partial^{\alpha}u\|^2+\|\partial^{\alpha}B\|^2   \right),\\
\label{GGG2.14}
&\mathcal{D}_1(t)=\|\nabla_x u\|_{H^2}^2+\|\nabla_x B\|_{H^2}^2.
\end{align}
Adding \eqref{HH2.7} to \eqref{HH2.9}, we have 
\begin{align}
\frac{\mathrm{d}}{\mathrm{d} t}\mathcal{E}_{1}(t)+\lambda_8\mathcal{D}_{1}(t)
\leq C\big((1+T)\delta_1+\delta_2\big)\mathcal{D}_{1}(t) 
\end{align}
for some $\lambda_8>0$.
Using the smallness of $\delta_1$ and $\delta_2$, we deduce that
\begin{align}\label{HHH2.14}
\frac{\mathrm{d}}{\mathrm{d} t}\mathcal{E}_{1}(t)+\lambda_9\mathcal{D}_{1}(t)
\leq 0 
\end{align}
for some $\lambda_9>0$.
%
%


In order to get the desired decay estimates to $(u,B)$, 
we consider the linearized equations  related to
the Cauchy problem\eqref{I3}-\eqref{I--3} as follows:
\begin{equation}\label{C1}\left\{
\begin{aligned}
& \partial_t u+\nabla_x P-\Delta_x u=0, \\
& \partial_t B-\Delta_x B=0, \\
& \nabla_x \cdot u=0, \quad \nabla_x \cdot B=0, \\
&  u(0, x)=u_0(x), \quad B(0, x)=B_0(x).
\end{aligned}\right.
\end{equation}
For brevity, 
we use $U(t) = (u(t), B(t))$ to represent the 
solution of the problem \eqref{C1} and $U_0 = (u_0, B_0)$.
Here, we denote $\|U\|_{L^2}=\|u\|_{L_x^2}+$ $\|B\|_{L_x^2}$, and $\|U\|_{L_1}=\|u\|_{L_x^1}+$ $\|B\|_{L_x^1}$.

By Duhamel's principle, we know the   solution to \eqref{C1} takes the form:
\begin{align}
U(t)=\Lambda(t) U_0,
\end{align}
where $\Lambda(t)$ is so-called the 
linear solution operator to \eqref{C1}.
Next, we give the main result on the  time decays of solution $U(t)$.
\begin{thm}\label{I6}
Let $1\leq q\leq2$. For any $\alpha$, $\alpha^{\prime}$ with 
$\alpha^{\prime}\leq\alpha$ and $m=|\alpha-\alpha^{\prime}|$,
\begin{equation}\label{D3}
\|\partial^{\alpha}\mathbb{A}(t)U_{0}\|_{L^2}
\leq C\left(1+t\right)^{-\frac{3}{2}\left(\frac{1}{q}-\frac{1}{2}\right)-\frac{m}{2}}
\left(\|\partial^{\alpha^{\prime}}U_{0}\|_{{L^q}}+\|\partial^{\alpha}U_{0}\|_{L^2}\right) 
\end{equation}
holds for all $t\geq 0$.
\end{thm}
\begin{proof}

Applying the Fourier transform to \eqref{C1} gives
\begin{equation}\label{G3.4}
\left\{
\begin{aligned}
& \partial_t \hat{u}+i k \hat{P}+|k|^2\hat{u}=0, \\
& \partial_t \hat{B}+|k|^2 \hat{B}=0, \\
&  k \cdot \hat{u}=0, \quad k \cdot \hat{B}=0 .
\end{aligned}
\right.
\end{equation}
Taking the product of \eqref{G3.4}$_1$ and \eqref{G3.4}$_2$ with $\overline{\hat{u}}$ and $\overline{\hat{B}}$ respectivly, we get
\begin{align}\label{G3.5}
\frac{1}{2} \frac{\mathrm{d}}{\mathrm{d} t}\left(|\hat{u}|^2+|\hat{B}|^2\right)+|k|^2|\hat{u}|^2+|k|^2|\hat{B}|^2+|k\cdot\hat{u}|^2+|k\cdot\hat{B}|^2 =0 .
\end{align}
We define
\begin{align}\label{G3.6}
\mathcal{E}\big(\hat{U}(t, k)\big)=|\hat{u}|^2+|\hat{B}|^2.
\end{align}
Substituting \eqref{G3.6} into \eqref{G3.5}, we get
\begin{align}
& \quad\frac{\mathrm{d}}{\mathrm{d} t}\mathcal{E}\big(\hat{U}(t, k)\big)+\frac{\lambda_{10}|k|^2}{1+|k|^2}\mathcal{E}\big(\hat{U}(t, k)\big) \leq 0 
\end{align}
for some $\lambda_{10}>0$, which  implies that
\begin{equation}
\mathcal{E}\big(\hat{U}(t, k)\big) \leq C \mathrm{e}^{\frac{-\lambda_{10}|k|^2 t}{1+|k|^2}} \mathcal{E}\big(\hat{U}(0, k)\big),
\end{equation}
where the Gronwall's inequality has been used.

As mentioned in \cite{DF-jmp-2010,Ksc-jjam-1984}, 
from the above estimate, we can easily get
the expected time-decay property \eqref{D3}.
For simplicity, we omit the details of the proof. 
\end{proof}

\medskip
Now we are  in a position   to prove the time-decay of Theorem \ref{I4}.
\begin{proof}[Proof of the decay part of Theorem \ref{I4}]
From \eqref{GGG2.13} and \eqref{GGG2.14}, we get
\begin{align}\label{G3.9}
\mathcal{E}_1(t)\leq C\big(\mathcal{D}_1(t)+\|u\|_{L^2}^2+\|B\|_{L^2}^2 \big). 
\end{align}
Adding $\|U\|_{L^2}^{2}$ 
to both sides of \eqref{H2.14} and using \eqref{G3.9} yield
\begin{align}\label{D11}
\frac{\rm d}{{\rm d}t}\mathcal{E}_{1}(t)+\lambda_{11}\mathcal{E}_{1}(t)\leq C\| U\|_{L^2}^{2} 
\end{align}
for some $\lambda_{11}>0$.

Therefore, we need to  give the estimate of $\|U\|_{L^2}^{2}$. 
By Duhamel’s principle, the solution $(u, B)$ to \eqref{I3}$_2$--\eqref{I--3} can be represented as
\begin{align}\label{gradU}
U(t)=\mathbb{A}(t)U_{0}+\int_{0}^{t}\mathbb{A}(t-s)(\mathfrak{G}_1,\mathfrak{G}_2) {\rm d}s:=\sum_{i=1}^{3}\mathfrak{J}_{i}(t),
\end{align}
where
\begin{align}
&\mathfrak{J}_{1}(t) =\mathbb{A}(t)U_{0}, \label{j-1} \\
&\mathfrak{J}_{2}(t) =\int_{0}^{t}\mathbb{A}(t-s)\left(\mathfrak{G}_1,0\right) {\rm d}s , \label{j-2} \\
&\mathfrak{J}_{3}(t) =\int_{0}^{t}\mathbb{A}(t-s)\left(0,\mathfrak{G}_2\right) {\rm d}s, \label{j-3}
\end{align}
with
\begin{align*}
&\mathfrak{G}_{1}=-u\cdot \nabla_x u+B\cdot\nabla_x B-\frac{1}{2}\nabla_x (|B|^2)+\int_{\mathbb{R}^3} 
\big((u-v)\times B\big)f\mathrm{d}v,\\
&\mathfrak{G}_2=B\cdot\nabla_x u-u\cdot\nabla_x B.
\end{align*}

Define ${\mathcal E}_{1,\infty}(t)=\sup\limits_{0\leq s\leq t}(1+s)^{\frac{3}{2}}{\mathcal E}_{1}(s)$. 
Then, by Theorem \ref{I6} and \eqref{HH2.1}, one has
\begin{align*}
\|\mathfrak{J}_{1}(t)\|_{L^2}&\leq C(1+t)^{-{\frac{3}{4}}}\|U_{0}\|_{{L^2}\bigcap {{L^1}}},\\
\|\mathfrak{J}_{2}(t)\|_{L^2}
&\leq C\int_{0}^{t}(1+t-s)^{-\frac{3}{4}}
\|G_{1}\|_{{L^2}\bigcap{L^1}} {\rm d}s \\
&\leq C\int_{0}^{t}(1+t-s)^{-\frac{3}{4}}\left(\|f\|_{H_{x,v}^{2}}{\mathcal E}_{1}^{\frac{1}{2}}(s)
+{\mathcal E}_{1}(s)
+\|f\|_{H_{x,v}^{2}}{\mathcal E}_{1}(s)\right) {\rm d}s \\
&\leq C\int_{0}^{t}(1+t-s)^{-\frac{3}{4}}(1+s)^{-\frac{3}{4}}\delta_1(s){\rm d}s
{\mathcal E}_{1,\infty}^{\frac{1}{2}}(t)\\
&\quad+C\int_{0}^{t}(1+t-s)^{-\frac{3}{4}}(1+s)^{-\frac{3}{2}} {\rm d}s {\mathcal E}_{1,\infty}(t)\\
&\leq C(1+t)^{-\frac{3}{4}}\left({{\mathcal E}_{1,\infty}^{\frac{1}{2}}(t)+\mathcal E}_{1,\infty}(t)\right),\\
\|\mathfrak{J}_{3}(t)\|_{L^2}
&\leq C\int_{0}^{t}(1+t-s)^{-\frac{3}{4}}\|G_{2}\|_{{L^2}\bigcap{L^1}} {\rm d}s \\
&\leq C\int_{0}^{t}(1+t-s)^{-\frac{3}{4}}\left({\mathcal E}_{1}(s)
\right) {\rm d}s \\
&\leq C(1+t)^{-\frac{3}{4}}{\mathcal E}_{1,\infty}(t).
\end{align*}
Thus, we   deduce that
\begin{align*}
\|U\|_{L^2}^{2}\leq\sum_{i=1}^{3}\|\mathfrak{J}_{i}(t)\|_{L^2}^{2}
\leq C(1+t)^{-\frac{3}{2}}\left(\|U_{0}\|_{{L^2}\bigcap{L^1}}^{2}+\mathcal{E}_{1,\infty}(t)
+\mathcal{E}_{1,\infty}^{2}(t)\right).
\end{align*}
Thanks to \eqref{D11}, we have
\begin{align*}
\frac{\rm d}{{\rm d}t}\mathcal{E}_{1}(t)+\lambda_{11}\mathcal{E}_{1}(t)
\leq C(1+t)^{-\frac{3}{2}}\left(\| U_{0}\|_{{L^2}\bigcap{L^1}}^{2}
+\mathcal{E}_{1,\infty}(t)+\mathcal{E}_{1,\infty}^{2}(t)\right).
\end{align*}
Using Gronwall’s inequality, we arrive at
\begin{align*}
\mathcal{E}_{1}(t)\leq e^{-\lambda_{11} t}\mathcal{E}_{1}(0)+C(1+t)^{-\frac{3}{2}}
\left(\| U_{0}\|_{{L^2}\bigcap{L^1}}^{2}+\mathcal{E}_{1,\infty}(t)
+\mathcal{E}_{1,\infty}^{2}(t)\right).
\end{align*}
Therefore,
\begin{align*}
\mathcal{E}_{1,\infty}(t)\leq C\left(\| U_{0}\|_{{H}^{2}\bigcap{L^1}}^{2}
+\mathcal{E}_{1,\infty}(t)+\mathcal{E}_{1,\infty}^{2}(t)\right).
\end{align*}
Since $\|U_0\|_{{H}^{2}\bigcap{L^1}}$ is small enough, we have
\begin{align*}
\mathcal{E}_{1,\infty}(t)\leq C\|U_{0}\|_{{H}^{2}\bigcap{L^1}}^{2} 
\end{align*}
for all $0\leq t\leq T$, which implies that
\begin{equation*}
\mathcal{E}_{1}(t)\leq C(1+t)^{-\frac{3}{2}}\|U_0\|_{{H}^{2}\bigcap{L^1}}^{2}.
\end{equation*}
Thus, we have
\begin{equation}\label{ubH2}
\|(u,B)\|_{{H}^{2}}\leq C(1+t)^{-\frac{3}{4}} 
\end{equation}
for all $0\leq t\leq T$.

Next, we begin to prove the time decay rate \eqref{G1.7}. We define high-order energy functional $\mathcal{H}_1(t)$ and high-order dissipation rate functional $\mathcal{M}_1(t)$  by 
\begin{align}\label{D12}
\mathcal{H}_1(t)&=\sum\limits_{1\leq|\alpha|\leq 2}
\big(\|\partial^{\alpha}u\|^{2}+\|\partial^{\alpha}B\|^{2}\big),\\
\mathcal{M}_1(t)&=\sum\limits_{1\leq|\alpha|\leq 2}
\big(\|\nabla_x\partial^{\alpha}u\|^{2}+\|\nabla_x\partial^{\alpha}B\|^{2}\big).
\end{align}
By using Lemma \ref{L1}, \eqref{HH2.1} and similar arguments to those in the proofs of 
Lemmas \ref{L2.3}--\ref{L2.4}, we get
\begin{align}\label{D14}
\frac{\rm d}{{\rm d}t}\mathcal{H}_{1}(t)+\lambda_{12}\mathcal{M}_{1}(t)\leq 0 
\end{align}
for some $\lambda_{12}>0$.

Adding the term $\|\nabla U\|_{L^2}^{2}$ to both sides of \eqref{D14} gives
\begin{align}\label{hha}
\frac{\rm d}{{\rm d}t}\mathcal{H}_{1}(t)+\lambda_{13}\mathcal{H}_{1}(t)\leq C\|\nabla U\|_{L^2}^{2} 
\end{align}
for some $\lambda_{13}>0$.

Applying  $\nabla$ to \eqref{gradU} gives
\begin{align}
  \nabla U(t)=\nabla \mathbb{A}(t)U_{0}+\int_{0}^{t}\nabla \mathbb{A}(t-s)(\mathfrak{G}_1,\mathfrak{G}_2) {\rm d}s:=\sum_{i=1}^{3}\nabla \mathfrak{J}_{i}(t),
\end{align}
where $ \mathfrak{J}_{i}(t) (i=1,2,3)$ are given in \eqref{j-1}--\eqref{j-3}.
Taking the same arguments as before, we get
\begin{align*}
\|\nabla \mathfrak{J}_{1}(t)\|_{L^2}&\leq C(1+t)^{-{\frac{5}{4}}}\|\nabla U_{0}\|_{{L^2}\bigcap {{L^1}}},\\
\|\nabla \mathfrak{J}_{2}(t)\|_{L^2}
&\leq C\int_{0}^{t}(1+t-s)^{-\frac{5}{4}}
\|\nabla G_{1}\|_{{L^2}\bigcap{L^1}} {\rm d}s \\
&\leq C\int_{0}^{t}(1+t-s)^{-\frac{5}{4}}\left(\|f\|_{H_{x,v}^{2}}{\mathcal H}_{1}^{\frac{1}{2}}(s)
+{\mathcal H}_{1}(s)
+\|f\|_{H_{x,v}^{2}}{\mathcal H}_{1}(s)\right) {\rm d}s \\
&\leq C\int_{0}^{t}(1+t-s)^{-\frac{5}{4}}(1+s)^{-\frac{5}{4}}{\rm d}s
{\mathcal H}_{1,\infty}^{\frac{1}{2}}(t)\\
&\quad+C\int_{0}^{t}(1+t-s)^{-\frac{5}{4}}(1+s)^{-\frac{5}{2}} {\rm d}s {\mathcal H}_{1,\infty}(t)\\
&\leq C(1+t)^{-\frac{5}{4}}\left({{\mathcal H}_{1,\infty}^{\frac{1}{2}}(t)+\mathcal H}_{1,\infty}(t)\right),\\
\|\nabla \mathfrak{J}_{3}(t)\|_{L^2}
&\leq C\int_{0}^{t}(1+t-s)^{-\frac{5}{4}}\|\nabla G_{2}\|_{{L^2}\bigcap{L^1}} {\rm d}s \\
&\leq C\int_{0}^{t}(1+t-s)^{-\frac{5}{4}}\left({\mathcal H}_{1}(s)
\right) {\rm d}s \\
&\leq C(1+t)^{-\frac{5}{4}}{\mathcal H}_{1,\infty}(t),
\end{align*}
where $\mathcal{H}_{1,\infty}(t)=\sup\limits_{0\leq s\leq t}(1+s)^{\frac{5}{2}}\mathcal{H}_{1}(s)$.
Therefore, we have
\begin{align}\label{D15}
\|\nabla U\|_{_{L^2}}^{2}\leq C(1+t)^{-\frac{5}{2}}
\left(\|U_0\|_{{H}^{2}\bigcap{L^1}}^{2}+\mathcal{H}_{1,\infty}(t)
+\mathcal{H}_{1,\infty}^{2}(t)\right). 
\end{align}

Combining \eqref{hha} and \eqref{D15} and using Gronwall's inequality, we arrive at
\begin{align*}
\mathcal{H}_{1}(t)\leq e^{-\lambda_{13} t}\mathcal{H}_{1}(0)
+C(1+t)^{-\frac{5}{2}}\left(\|U_0\|_{{H}^{2}\bigcap{L^1}}^{2}
+\mathcal{H}_{1,\infty}(t)+\mathcal{H}_{1,\infty}^{2}(t)\right),
\end{align*}
and hence
\begin{equation*}
\mathcal{H}_{1,\infty}(t)\leq C\left(\|U_0\|_{{H}^{2}\bigcap{L^1}}^{2}
+\mathcal{H}_{1,\infty}(t)+\mathcal{H}_{1,\infty}^{2}(t)\right).
\end{equation*}
Due to $\|U_0\|_{{H}^{2}\bigcap{L^1}}^{2}$ is small enough, which means
\begin{equation*}
\mathcal{H}_{1,\infty}(t)\leq C\|U_0\|_{{H}^{2}\bigcap{L^1}}^{2}
\end{equation*}
holds for all $0\leq t\leq T$. By the definition of $\mathcal{H}_{1,\infty}(t)$, we have
\begin{equation*}
\mathcal{H}_{1}(t)\leq C(1+t)^{-\frac{5}{2}}\|U_0\|_{{H}^{2}\bigcap{L^1}}^{2},
\end{equation*}
which implies that
\begin{equation*}
\|\nabla(u,B)\|_{{H}^{1}}\leq C(1+t)^{-\frac{5}{4}} 
\end{equation*}
for all $0\leq t\leq T$. Hence we completes the proof of \eqref{G1.7}.

Noticing that by applying the Fourier transform method, we know that 
the decay rate of $\|(u,B)\|_{{H}^{2}}$ is $(1+t)^{-\frac{3}{4}}$. However, through \eqref{G1.6}, we   deduce that
\begin{align*}
\|(u,B)\|_{{H}^{2}}\leq C\|(u_0,B_0)\|_{{H}^{2}}\leq C(1+T)^{-1}\leq C(1+t)^{-1}\leq C(1+t)^{-\frac{3}{4}}   
\end{align*}
for all $0\leq t\leq T$.
Hence  we finally complete the proof of Theorem \ref{I4}. 
\end{proof}

\smallskip
\section{The Vlasov-MHD system  near the equilibrium state $(0,0,0)$ in $\mathbb{T}^3$}

In this section, we study the Vlasov-MHD system \eqref{I3}  near the state $(0,0,0)$  in $\mathbb{T}^3$.
We shall establish the global existence of strong solutions when the initial data is sufficient small.
We also prove the exponential decay of $u $ and $B$. 

\begin{proof}[Proof of Theorem \ref{I5}] 
Due to the boundedness of the periodic domain,  we only need to  assume that $f_0(x,v)$ has a compact 
support  in the velocity $v$.  The existence and uniqueness of strong solutions to \eqref{I3} in the  periodic domain $\mathbb{T}^3$ case are essential to those in the  whole space $\mathbb{R}^3$ case and we omit them for brevity.  
 Below we focus on the uniform a priori estimates to the strong solutions to obtain the exponential decay of $u $ and $B$. 

First, we   obtain  the following conservation laws by a direct calculation:
\begin{equation}\nonumber
\left\{
\begin{aligned}
&\frac{\mathrm{d}}{\mathrm{d} t} \int_{\mathbb{R}^3\times\mathbb{T}^3} f \mathrm{d} v \mathrm{d} x=0,\\
& \frac{\mathrm{d}}{\mathrm{d} t}\left(\int_{\mathbb{T}^3}  u \mathrm{d} x+\int_{\mathbb{R}^3\times\mathbb{T}^3} v f \mathrm{d} v \mathrm{d} x\right)=0, \\
&\frac{\mathrm{d}}{\mathrm{d} t} \int_{\mathbb{T}^3} B \mathrm{d} x=0.
\end{aligned}\right.
\end{equation}
Under the assumption of Theorem \ref{I5}, we have
\begin{equation} \label{G4.1}
\left\{
\begin{aligned}
& \int_{\mathbb{R}^3\times\mathbb{T}^3} f \mathrm{d} v \mathrm{d} x=c_1, \\
& \int_{\mathbb{T}^3}  u \mathrm{d} x+\int_{\mathbb{R}^3\times\mathbb{T}^3} v f \mathrm{d} v \mathrm{d} x=c_2,  \\  
& \int_{\mathbb{T}^3} B \mathrm{d} x=0 
\end{aligned}\right.
\end{equation}
 for all $0 \leq t\leq T$,
where $c_1$ and $c_2$ are constants  satisfying $0<c_1,c_2\ll 1$. 

According to the Poincaré's inequality and  \eqref{G4.1}, we obtain that
\begin{align}\label{G4.2}
\|B\|_{L^2} & \leq C\|\nabla_x B\|_{L^2},\\
\|u\|_{L^2} & \leq C\Big\|\nabla_x u+\int_{\mathbb{R}^3}\bigg(v-\frac{c_2}{c_1}\bigg)\nabla_x f\mathrm{d}v\Big\|_{L^2}
+C\Big\|\int_{\mathbb{R}^3}\bigg(v-\frac{c_2}{c_1}\bigg)\nabla_x f\mathrm{d}v\Big\|_{L^2} \nonumber\\
&\leq C\|\nabla_x u\|_{L^2}+C\|f\|_{H_{x,v}^2}\nonumber \\
\label{G4.3}
&\leq C \|\nabla_x u\|_{L^2}+C\epsilon_3^{\frac{1}{2}} .
\end{align}
Combining \eqref{G4.2} with \eqref{G4.3}, one has 
\begin{align}\label{G4.5}
\|u\|_{L^2}^2+\|B\|_{L^2}^2\leq C\big(\|\nabla_x u\|_{L^2}^2+ \|\nabla_x B\|_{L^2}^2 \big)  +C\epsilon_3. 
\end{align}

We define the energy functionals $\mathcal{E}_1(t)$ and the corresponding dissipation rate functional $\mathcal{D}_1(t)$  to the system \eqref{I3} as same as the $\mathbb{R}^3$ case given in \eqref{GGG2.13}
and \eqref{GGG2.14}, respectively.  Similar to the previous arguments on the whole space case, it follows that
\begin{align}\label{G4.7}
\frac{\rm d}{{\rm d}t}\mathcal{E}_{1}(t)+\lambda_{14}\mathcal{D}_{1}(t)\leq 0 
\end{align}
for some $\lambda_{14}>0$.

Thus, we combine \eqref{G4.7} and \eqref{G4.5} to get 
\begin{align}\label{G4.6}
 \frac{\rm d}{{\rm d}t}\mathcal{E}_{1}(t)+\lambda_{15}\mathcal{E}_{1}(t)\leq C\epsilon_3    
\end{align}
for some $\lambda_{15}>0$.
Next, we define
\begin{align*}
\mathcal{Q}(t):=e^{\lambda_{15} t}\mathcal{E}_1(t).
\end{align*}
From \eqref{G4.6}, we   obtain that
\begin{align}\label{K4.7}
\frac{\mathrm{d}}{\mathrm{d} t}\mathcal{Q}(t)\leq e^{\lambda_{15} t} \left(\frac{\mathrm{d}}{\mathrm{d} t}\mathcal{E}_1(t)+\lambda_{15}\mathcal{E}_1(t)\right)\leq C e^{\lambda_{15} t}\epsilon_3.
\end{align}
Thus,
\begin{align*}
\mathcal{Q}(t)\leq \mathcal{Q}(0)+C\int_0^Te^{\lambda_{15} t}\epsilon_3     \mathrm{d}t. 
\end{align*}
Namely, we get
\begin{align*}
\mathcal{E}_1(t)&\leq e^{-\lambda_{15} t}\left(\mathcal{E}_1(0)+CTe^{\lambda_{15} T}\epsilon_3                 \right) \leq Ce^{-\lambda_{15} t}\Big(\epsilon_3^{\frac{1}{2}}+\epsilon_4\Big) 
 \leq Ce^{-\lambda_{15} t} ,
\end{align*}
where we have used the smallness of $\epsilon_3$.
Thus, we obtain 
\begin{align*}
& \|( u, B)\|_{H^2} \leq C e^{-\lambda_{16} t}    
\end{align*}
for some $\lambda_{16}>0$ and all $0\leq t\leq T$, and then complete the proof of Theorem \ref{I5}.
\end{proof} 

\section{The Vlasov-MHD system  near the equilibrium state  $(M,0,0)$ in $\mathbb{R}^3$ and $\mathbb{T}^3$}

In this section, we study the Vlasov-MHD system \eqref{I3} near the equilibrium state  $(M,0,0)$ 
in $\mathbb{R}^3$ and $\mathbb{T}^3$. 
Compared to the vacuum case   $(\bar f,\bar u,\bar B)=(0,0,0)$,
  for the case $( \tilde f,\tilde u,\tilde B)=(M,0,0)$,
   there are two  additional terms in the system \eqref{I8}, i.e.  $(u\times B)\cdot v M$ in the equation of $g$ 
    and $ u\times B$ in the equation of $u$, which bring  some difficulties in handling  the  upper bound of $\|g\|_{H^2_{x,v}}$. We need to assume that the initial values of $u$ and $B$ are controlled by an order of $e^{-\mathcal{O}(1)T}$. 

\subsection{The whole space $\mathbb{R}^3$ case}   In this case, due to the algebraic decay obtained by performing Fourier transform on the whole space
     are slower than the exponential decay given via the smallness of initial
     data, we don't need to consider the large time behavior here.
 
\begin{proof}[Proof of Theorem \ref{T5.2}]
The whole proof is   divided into two steps: (1) establish  the a priori estimates  of strong solutions;
(2) prove the existence and  uniqueness of strong solutions.

%
{\it Step 1: A priori estimates.} \ \
We shall establish the uniform  a priori estimates of $g$, $u$ and $B$. We first suppose that $(g,u,B)$ is the strong solution to the Cauchy problem \eqref{I8}-\eqref{I--8} on $0\leq t\leq T$, and it satisfies 
\begin{align}\label{N5.3}
\sup_{ {0\leq t\leq T}}\Big(\|g\|_{H_{x,v}^2}+\|(u,B)\|_{H^{2}}\Big)\leq& \delta_3 
\end{align}
with $\delta_3$ taking an order of $e^{-\mathcal{O}(1)T}$.

Let $\tilde\Sigma_x (t)$ and $\tilde\Sigma _v(t)$ be the $x$-- and $v$-- projections of supp$g(t,\cdot,\cdot)$ respectively,
where
\begin{align*}
\tilde{\Sigma}_x (t):=&\{x\in\mathbb{R}^3_x:\:\exists\,(x,v)\in\mathbb{R}^3_x\times\mathbb{R}^3_v,\  \text{such that}\  g(t, x,v)\neq0\},\\
\tilde{\Sigma}_v (t):=&\{v\in\mathbb{R}_v^3:\:\exists\,(x,v)\in\mathbb{R}^3_x\times\mathbb{R}^3_v,\  \text{such that}\  g(t ,x,v)\neq0\}.
\end{align*}
In addition,
we define $\tilde{\mathcal{R}}_x ( t) $ and $\tilde{\mathcal{R}}_v ( t)$  as follows:
\begin{align*}
\tilde{\mathcal{R}}_{x} (t):=\sup_{x\in\tilde{ \Sigma}_{x} (t)}|x|,
\quad \tilde{\mathcal{R}}_{v}(t):=\sup_{v\in \tilde{\Sigma}_{v}(t)}|v| .
\end{align*}

\begin{lem}\label{L5.2}
For any $T>0$, suppose that $g_0(x,v)$ has a compact 
support in $x$ and $v$, it holds that
\begin{align}
\tilde{\mathcal{R}}_{x}(t)&\leq \tilde{\mathcal{R}}_{x}(0)+C\left(\tilde{\mathcal{R}}_{v}(0)+\delta_3 T\right)T,\\
\tilde{\mathcal{R}}_{v}(t)&\leq C\big(\tilde{\mathcal{R}}_{v}(0)+\delta_3 T\big) 
\end{align}
for all $0\leq t\leq T$.
\end{lem}

\begin{proof}
Taking  the process of Lemma \ref{N2.3},
we  achieve the Lemma \ref{L5.2} immediately.
\end{proof}

Thanks to the   smallness of $\delta_3$, we  also obtain  
\begin{align}\label{N5.6}
\tilde{\mathcal{R}}_{x}^\infty \approx C(1+T),\ \ \ \ \tilde{\mathcal{R}}_{v}^\infty \approx C.
\end{align} 

For the simplicity of symbols, we set 
\begin{align*}
\tilde{\mathcal{R}}_{x}^{\infty}&:=\sup_{0\leq t\leq T_1}\tilde{\mathcal{R}}_{x}(t),\quad \tilde{\mathcal{B}}_{x}=\tilde{\mathcal{B}}_{\tilde{\mathcal{R}}_{x}^{\infty}},
\\ \tilde{\mathcal{R}}_{v}^{\infty}&:=\sup_{0\leq t\leq T_1}R_{v}(t),\quad \tilde{\mathcal{B}}_{v}=\tilde{\mathcal{B}}_{\tilde{\mathcal{R}}_{v}^{\infty}},
\end{align*}
where $\tilde{\mathcal{B}}_{\tilde{\mathcal{R}}_{x}^{\infty}}$ and $\tilde{\mathcal{B}}_{\tilde{\mathcal{R}}_{v}^{\infty}}$ denote the ball of radius $\tilde{\mathcal{R}}_{x}^{\infty}$ and $\tilde{\mathcal{R}}_{v}^{\infty}$, respectively.

Below we give some useful lemmas.
\begin{lem}\label{L5.3}
For the chosen $T$ stated in Lemma \ref{N5.3}, there exists a positive constant $\lambda_{17}>0$, such that
\begin{align}
&\frac{1}{2} \frac{\mathrm{d}}{\mathrm{d} t}\left(\|u\|^2+\|B\|^2\right)+\lambda_{17}\left(\|\nabla_x u\|^2+\|\nabla_x B\|^2\right)\nonumber \\ \label{N5.9}
&\qquad \quad \leq C\delta_3(1+T)\left(\|\nabla_x u\|^2+\|\nabla_x B\|^2\right)
\end{align}
holds for all $0\leq t\leq T$.
\end{lem}
\begin{proof}
Multiplying the equations \eqref{I8}$_2$--\eqref{I8}$_3$ by $u$ and $B$ respectively, 
and then taking integration and summation, we obtain
\begin{align}
&\frac{1}{2} \frac{\mathrm{d}}{\mathrm{d} t}\left(\|u\|^2+\|B\|^2\right)+\|\nabla_x (u,B)\|^2\nonumber\\
=&-\int_{\mathbb{R}^3} u\cdot\nabla_x u\cdot u\mathrm{d}x-\int_{\mathbb{R}^3} u\cdot\nabla_x B\cdot B\mathrm{d}x+\int_{\mathbb{R}^3} B\cdot\nabla_x u\cdot B\mathrm{d}x\nonumber\\
&+\int_{\mathbb{R}^3} B\cdot\nabla_x B\cdot u\mathrm{d}x-\frac{1}{2}\int_{\mathbb{R}^3} \nabla_x (|B|^2)\cdot u\mathrm{d}x\nonumber\\
&+\int_{\mathbb{R}^3\times\mathbb{R}^3} 
\big(u\cdot(u-v)\times B\big)g\mathrm{d}v\mathrm{d}x+\int_{\mathbb{R}^3} u\cdot(u \times B)\mathrm{d}x\nonumber\\
:=&\sum_{j=15}^{21}\mathcal{I}_{j}.
\end{align}
For the terms $\mathcal{I}_j$ ($j={15},\dots,{19}$), 
using the inequalities of Lemma \ref{L1}, we have
\begin{align*}
 \mathcal{I}_{15}\leq\,  &C\|u\|_{L^3}\|\nabla_x u\|\|u\|_{L^6}\ \leq C\|u\|_{H^1}\|\nabla_x u\|^2\leq C \delta_3\|\nabla_x u\|^2,\\
 \mathcal{I}_{16}\leq\,  & C\|u\|_{L^3}\|\nabla_x B\|\|B\|_{L^6}\leq C\|u\|_{H^1}\|\nabla_x B\|^2\leq C \delta_3\|\nabla_x B\|^2,\\
 \mathcal{I}_{17}\leq \,  &C\|B\|_{L^3}\|\nabla_x u\|\|B\|_{L^6}\leq C\|B\|_{H^1}\|\nabla_x (u,B)\|^2\leq C \delta_3\|\nabla_x (u,B)\|^2,\\
 \mathcal{I}_{18}\leq \,  &C\|u\|_{L^3}\|\nabla_x B\|\|B\|_{L^6}\leq C\|u\|_{H^1}\|\nabla_x B\|^2\leq C \delta_3\|\nabla_x B\|^2,\\
 \mathcal{I}_{19}\leq \,  &C\|u\|_{L^3}\|\nabla_x B\|\|B\|_{L^6}\leq C\|u\|_{H^1}\|\nabla_x B\|^2\leq C \delta_3\|\nabla_x B\|^2.
\end{align*}
For the term $\mathcal{I}_{20}$, using H\"{o}lder’s inequality and \eqref{N5.6}, we have
\begin{align*}
\mathcal{I}_{20}\,&\leq  C\int_{\mathbb{R}^3\times \tilde{\mathcal{B}}_v}|u|^2|B||g|_2\mathrm{d}v\mathrm{d}x+C\int_{\tilde{\mathcal{B}}_x\times \tilde{\mathcal{B}}_v}|u||B||g|_2\mathrm{d}v\mathrm{d}x\ \\
&\leq C\|u\|_{L^6}^2\|B\|_{L^6}\|g\|_{H_{x,v}^2}
+C\|u\|_{L^6}\|B\|_{L^6}\|1{\chi({\tilde{\mathcal{B}}_x})}\|_{L^6(\tilde{\mathcal{B}}_{x})}\|g\|_{H_{x,v}^2}    \\
&\leq C\|g\|_{H_{x,v}^2}\|u\|_{H^2}\left(\|\nabla_x u\|^2+\|\nabla_x B\|^2 \right)+C(1+T)\|g\|_{H_{x,v}^2}\left(\|\nabla_x u\|^2+\|\nabla_x B\|^2 \right)\\
&\leq C\delta_3\|\nabla_x u\|^2+C\delta_3(1+T) \big(\|\nabla_x u\|^2+\|\nabla_x B\|^2 \big).
\end{align*}
For the last term $\mathcal{I}_{21}$, using the identity $u\cdot(u\times B)=-B\cdot(u\times u)=0$, 
we get $\mathcal{I}_{21}=0.$

Collecting the above estimates together, we  easily get the estimate \eqref{N5.9}. We now complete the proof of Lemma \ref{L5.3}.
\end{proof}

\begin{lem}\label{L5.4}
For the chosen $T$ stated in Lemma \ref{N5.3}, there exists a positive constant $\lambda_{18}>0$, such that
\begin{align} 
&  \frac{1}{2} \frac{\mathrm{d}}{\mathrm{d} t} \sum_{1 \leq|\alpha| \leq 2}\left(\|\partial^\alpha u\|^2+\|\partial^\alpha B\|^2\right)+\lambda_{18} \sum_{1 \leq|\alpha| \leq 2}\Big(\|\nabla_x \partial^\alpha u\|^2+\|\nabla_x \partial^\alpha B\|^2 \Big)\nonumber\\\label{N5.11}
&\quad\quad\quad\leq C\delta_3\left(\|\nabla_x u\|_{H^2}^2+\|\nabla_x B\|_{H^2}^2\right) 
\end{align}
holds for all $0\leq t\leq T$.
\end{lem}
\begin{proof}
Applying $\partial^{\alpha}$ with $1 \leq|\alpha| \leq 2$ 
to the equations \eqref{I8}$_2$--\eqref{I8}$_3$  
and multiplying the results by $\partial^\alpha u$ and $ \partial^\alpha B$ respectively, then taking integration and summing them together, we  obtain
\begin{align}
& \frac{1}{2} \frac{\mathrm{d}}{\mathrm{d} t} \left(\left\|\partial^\alpha u\right\|^2+\left\|\partial^\alpha B\right\|^2\right)+\|\nabla_x \partial^\alpha u\|^2+\|\nabla_x \partial^\alpha B\|^2 \nonumber\\
=&-\int_{\mathbb{R}^3}\partial^{\alpha}(u\cdot \nabla_x B)\cdot\partial^{\alpha}B\mathrm{d}x+\int_{\mathbb{R}^3}\partial^{\alpha}(B\cdot \nabla_x B)\cdot\partial^{\alpha}u\mathrm{d}x\nonumber\\
&+\int_{\mathbb{R}^3}\partial^{\alpha}(B\cdot \nabla_x u)\cdot\partial^{\alpha}B\mathrm{d}x-\frac{1}{2}\int_{\mathbb{R}^3}\partial^{\alpha}\nabla_x (|B|^2)\cdot\partial^{\alpha}u\mathrm{d}x-\int_{\mathbb{R}^3}\partial^{\alpha}(u\cdot \nabla_x u)\cdot\partial^{\alpha}u\mathrm{d}x\nonumber\\
&+\int_{\mathbb{R}^3\times\mathbb{R}^3} 
\partial^{\alpha}u\cdot\partial^{\alpha}\Big\{\big((u-v)\times B\big)g\Big\}\mathrm{d}v\mathrm{d}x+\int_{\mathbb{R}^3}\partial^{\alpha}(u\times B)\cdot\partial^{\alpha}u\mathrm{d}x
\nonumber\\
:=&\sum_{j=22}^{28}\mathcal{I}_{j} .
\end{align}
With the help of Lemma \ref{L1}, we have
\begin{align*}
\mathcal{I}_{22}\,&\leq C\|\partial^{\alpha}(u\cdot\nabla_x B)\|\|\partial^{\alpha} B\|\\
&\leq C\|B\|_{H^2}\left(\|\nabla_x u\|_{H^2}^2+\|\nabla_x B\|_{H^2}^2\right)\\
&\leq C\delta_3\left(\|\nabla_x u\|_{H^2}^2+\|\nabla_x B\|_{H^2}^2\right).
\end{align*}
Similarly, we   achieve
\begin{align*}
\mathcal{I}_{23}&\leq C\|u\|_{H^2}\|\nabla_x B\|_{H^2}^2\leq C\delta_3\|\nabla_x B\|_{H^2}^2,\\
\mathcal{I}_{24}&\leq C\|B\|_{H^2}\left(\|\nabla_x B\|_{H^2}^2+\|\nabla_x u\|_{H^2}^2     \right)\leq C\delta_3\left(\|\nabla_x B\|_{H^2}^2+\|\nabla_x u\|_{H^2}^2     \right),\\
\mathcal{I}_{25}&\leq C\|u\|_{H^2}\|\nabla_x B\|_{H^2}^2\leq C\delta_3\|\nabla_x B\|_{H^2}^2,\\
\mathcal{I}_{26}&\leq C\|u\|_{H^2}\|\nabla_x u\|_{H^2}^2\leq C\delta_3\|\nabla_x u\|_{H^2}^2,\\
\mathcal{I}_{27}&\leq C\|u\|_{H^2}\left(\|\nabla_x B\|_{H^2}^2+\|\nabla_x u\|_{H^2}^2     \right)\leq C\delta_3\left(\|\nabla_x B\|_{H^2}^2+\|\nabla_x u\|_{H^2}^2     \right).
\end{align*}
For the remaining term $\mathcal{I}_{27}$, 
using H\"{o}lder’s inequality and Lemma \ref{L1}, 
we get
\begin{align*}
\mathcal{I}_{27}\,&\leq C\int_{\mathbb{R}^3\times\mathbb{R}^3} 
|\partial^{\alpha}u||\partial^{\alpha}\big(( u\times B)g\big)|\mathrm{d}v\mathrm{d}x+C\int_{\mathbb{R}^3\times\mathbb{R}^3} 
|\partial^{\alpha}u||\partial^{\alpha}\big(( v\times B)g\big)|\mathrm{d}v\mathrm{d}x\\
&\leq C\|\nabla_x (u\times B)\|_{H^1}\|g\|_{H_{x,v}^2}\|\nabla_x u\|_{H^2}+C\|\nabla_x B\|_{H^1}\|g\|_{H_{x,v}^2}\|\nabla_x u\|_{H^2}\\
&\leq C\|B\|_{H^2}\|g\|_{H_{x,v}^2}\|\nabla_x u\|_{H^2}^2+C\|g\|_{H_{x,v}^2}\left(\|\nabla_x u\|_{H^2}^2+\|\nabla_x B\|_{H^2}^2\right)\\
&\leq C\delta_3\left(\|\nabla_x u\|_{H^2}^2+\|\nabla_x B\|_{H^2}^2\right).
\end{align*}
Combining the above estimates, and then taking summation over $1 \leq|\alpha| \leq 2$, we obtain the estimate \eqref{N5.11}, and hence  finish the proof of Lemma \ref{L5.4}.
\end{proof}

Next, we define the temporal energy functional 
$\mathcal{E}_2(t)$ and the corresponding dissipation rate
functional $\mathcal{D}_2(t)$ to the system \eqref{I8} as follows:
\begin{align}\label{N5.13}
&\mathcal{E}_2(t)=\|u\|^2+\|B\|^2+\sum_{1\leq |\alpha|\leq 2}\left(\|\partial^{\alpha}u\|^2+\|\partial^{\alpha}B\|^2   \right),\\
\label{N5.14}
&\mathcal{D}_2(t)=\|\nabla_x u\|_{H^2}^2+\|\nabla_x B\|_{H^2}^2.
\end{align}
Adding \eqref{N5.9} and \eqref{N5.11}, we obtain 
\begin{align}
\frac{\mathrm{d}}{\mathrm{d} t}\mathcal{E}_{2}(t)+\lambda_{19}\mathcal{D}_{2}(t)
\leq C\delta_3\mathcal{D}_{2}(t)
\end{align}
for some $\lambda_{19}>0$.
Using the smallness of $\delta_3$, we  infer that
\begin{align}\label{N5.16}
\frac{\mathrm{d}}{\mathrm{d} t}\mathcal{E}_{2}(t)+\lambda_{20}\mathcal{D}_{2}(t)
\leq 0
\end{align}
for some $\lambda_{20}>0$.
Integrating \eqref{N5.16} over $(0,t)$ yields
\begin{align}\label{N5.17}
\mathcal{E}_{2}(t)+\lambda_{20}\int_0^t \mathcal{D}_{2}(s)\mathrm{d} s\leq \mathcal{E}_{2}(0)    
\end{align}
for all $0\leq t \leq T$.

Therefore, we have 
\begin{align}\label{G5.17}
\int_0^t \big(\|\nabla_x u(s)\|_{H^2}^2+\|\nabla_x B(s)\|_{H^2}^2\big)\mathrm{d} s\leq \mathcal{E}_{2}(0) 
\end{align}
for all $0\leq t \leq T$.

\medskip 
With the above preparations,
below we will give the estimate of $g$.

\begin{lem}\label{L5.7}
For the chosen $T$ stated in Lemma \ref{N5.3}, the following inequality  
\begin{align}\label{G5.18}
&\frac{1}{2} \frac{\mathrm{d}}{\mathrm{d} t}\|g\|_{H_{x,v}^2}
\leq C\big(1+\|\nabla_x u\|_{H^2}^2+\|\nabla_x B\|_{H^2}^2\big)\|g\|_{H_{x,v}^2}+C \|u\|_{H^2}\|B\|_{H^2}
\end{align}
holds for all $0\leq t\leq T$.
\end{lem}
\begin{proof}
Firstly, we consider the zeroth-order estimate of $g$. 
Multiplying \eqref{I8}$_1$ by $g$ and integrating over $\mathbb{R}^3\times\mathbb{R}^3$, one gets
\begin{align*}
\frac{1}{2} \frac{\mathrm{d}}{\mathrm{d} t}\|g\|^2= -\int_{\mathbb{R}^3\times\mathbb{R}^3}(u\times B)\cdot vMg\mathrm{d}v\mathrm{d}x 
\leq   C \|u\|_{H^2}\|B\|_{H^2}\|g\|_{H_{x,v}^2}.
\end{align*}

Secondly, we study the pure spatial derivatives of $g$. 
Applying $\partial^{\alpha}$ with $1 \leq|\alpha| \leq 2$ 
to the equation \eqref{I8}$_1$, we get
\begin{align}\label{G5.19}
&\partial^\alpha \partial_t g+v\cdot\nabla_x\partial^\alpha g+\sum_{|\alpha^\prime|\leq |\alpha|}C_{\alpha}^{\alpha^\prime}v\times
\partial^{\alpha^\prime}B\cdot\nabla_v\partial^{\alpha-\alpha^\prime}g\nonumber \\
&\qquad -\sum_{|\alpha^\prime|\leq |\alpha|}C_{\alpha}^{\alpha^\prime}\partial^{\alpha^\prime}(u\times B)\cdot\nabla_v\partial^{\alpha-\alpha^\prime}g+\partial^{\alpha}(u\times B)\cdot vM=0.
\end{align}
Then, multiplying \eqref{G5.19}
by $\partial^{\alpha}g$, 
and taking integration, we get
\begin{align*}
\frac{1}{2}\frac{\mathrm{d}}{\mathrm{d}t}\|\partial^\alpha g\|^2=&\sum_{1\leq|\alpha^\prime|\leq |\alpha|}C_{\alpha}^{\alpha^\prime}\int_{\mathbb{R}^3\times\mathbb{R}^3}
\big(v\times\partial^{\alpha^\prime}B\cdot\nabla_v\partial^{\alpha-\alpha^\prime}g \big)\partial^\alpha g\mathrm{d}v\mathrm{d}x \\
&-\sum_{1\leq|\alpha^\prime|\leq |\alpha|}C_{\alpha}^{\alpha^\prime}\int_{\mathbb{R}^3\times\mathbb{R}^3}\big(\partial^{\alpha^\prime}(u\times B)\cdot\nabla_v\partial^{\alpha-\alpha^\prime}g\big)\partial^\alpha g\mathrm{d}v\mathrm{d}x\\
&-\int_{\mathbb{R}^3\times\mathbb{R}^3}\partial^{\alpha}(u\times B)\cdot vM
\partial^{\alpha}g\mathrm{d}v\mathrm{d}x   \\
\leq\,& C\int_{\mathbb{R}^3\times B_{v}}\big(1+|\partial^{\alpha^\prime}B|^2\big)|\nabla_v\partial^{\alpha-\alpha^\prime}f||\partial^\alpha g|\mathrm{d}v\mathrm{d}x\\
&+C\int_{\mathbb{R}^3\times \mathbb{R}^3}|\partial^{\alpha^\prime}(u\times B)||\nabla_v\partial^{\alpha-\alpha^\prime}g||\partial^\alpha g|\mathrm{d}v\mathrm{d}x+C\|\partial^{\alpha}(u\times B)\| \|\partial^{\alpha}g\|\\
\leq\,& C\big(\|\nabla_x B\|_{H^2}\| u\|_{H^2} +\|\nabla_x u\|_{H^2}\| B\|_{H^2}  \big)\|g\|_{H_{x,v}^2}^2\\
&+C \|g\|_{H_{x,v}^2}^2+C\|\nabla_x B\|_{H^2}^2\|g\|_{H_{x,v}^2}^2+C\|u\|_{H^2}\|B\|_{H^2}\|g\|_{H_{x,v}^2}\\
\leq\,& C\big(1+\|\nabla_x u\|_{H^2}^2+\|\nabla_x B\|_{H^2}^2\big)\|g\|_{H_{x,v}^2}^2+C\|u\|_{H^2}\|B\|_{H^2}\|g\|_{H_{x,v}^2},
\end{align*}
where we have employed 
H\"{o}lder’s inequality and Lemma \ref{L1}.

Finally, we estimate the space-velocity-mixed 
derivatives of $g$. 
Select $\alpha,\beta$ with $|\alpha|+|\beta|\leq 2$.
Multiplying \eqref{I8}$_1$ by $\partial_{\beta}^\alpha$, one has
\begin{align*}
&\partial_{\beta}^\alpha \partial_t g+\sum_{\beta^\prime\neq 0}C_{\beta}^{\beta^\prime}
\partial_{\beta^\prime}v\cdot\nabla_v\partial_{\beta-\beta^\prime}^{\alpha}g-\sum_{\alpha^\prime\neq 0}C_{\alpha}^{\alpha^\prime}\partial^{\alpha^\prime}(u\times B)\cdot\nabla_v\partial_{\beta}^{\alpha-\alpha^\prime}g\\
&\qquad\quad =(u\times B)\cdot \nabla_v \partial_{\beta}^\alpha g-\sum_{\beta^\prime\neq 0}C_{\alpha}^{\alpha^\prime}C_{\beta}^{\beta^\prime}\partial_{\beta^\prime}v\times\partial^{\alpha^\prime} B\cdot\nabla_v\partial_{\beta-\beta^{\prime}}^{\alpha-\alpha^\prime}g+\partial^{\alpha}(u\times B)\cdot \partial_{\beta}(vM).  
\end{align*}
Taking integration, one obtains
\begin{align*}
\frac{1}{2}\frac{\mathrm{d}}{\mathrm{d}t}\|\partial_{\beta}^\alpha g\|^2=&-\sum_{\beta^\prime\neq 0}C_{\beta}^{\beta^\prime}\int_{\mathbb{R}^3\times\mathbb{R}^3}
(\partial_{\beta^\prime}v\cdot\nabla_v\partial_{\beta-\beta^\prime}^{\alpha}g)\partial_{\beta}^\alpha g\mathrm{d}v\mathrm{d}x\\
&+\sum_{\alpha^\prime\neq 0}C_{\alpha}^{\alpha^\prime}\int_{\mathbb{R}^3\times\mathbb{R}^3}\big(\partial^{\alpha^\prime}(u\times B)\cdot\nabla_v\partial_{\beta}^{\alpha-\alpha^\prime}g\big)\partial_{\beta}^\alpha g\mathrm{d}v\mathrm{d}x\\
&-\sum_{\beta^\prime\neq 0}C_{\alpha}^{\alpha^\prime}C_{\beta}^{\beta^\prime}
\int_{\mathbb{R}^3\times\mathbb{R}^3}\big(\partial_{\beta^\prime}v\times\partial^{\alpha^\prime} B\cdot\nabla_v\partial_{\beta-\beta^{\prime}}^{\alpha-\alpha^\prime}g\big)\partial_{\beta}^\alpha g\mathrm{d}v\mathrm{d}x\\
&-\int_{\mathbb{R}^3\times \mathbb{R}^3}\partial^{\alpha}(u\times B)\cdot\partial_{\beta}(vM)\partial_{\beta}^{\alpha}g\mathrm{d}v\mathrm{d}x\\
\leq\,&C\|g\|_{H_{x,v}^2}^2+C\big(\|\nabla_x B\|_{H^2}\| u\|_{H^2} +\|\nabla_x u\|_{H^2}\| B\|_{H^2}  \big)\|g\|_{H_{x,v}^2}^2\\
&+C \|g\|_{H_{x,v}^2}^2+C\|\nabla_x B\|_{H^2}^2\|g\|_{H_{x,v}^2}^2+C\|u\|_{H^2}\|B\|_{H^2}\|g\|_{H_{x,v}^2}\\
\leq\,& C\big(1+\|\nabla_x u\|_{H^2}^2+\|\nabla_x B\|_{H^2}^2\big)\|g\|_{H_{x,v}^2}^2+C\|u\|_{H^2}\|B\|_{H^2}\|g\|_{H_{x,v}^2}.
\end{align*}
Combining the above estimates and taking summation,
we   directly obtain the desired estimate \eqref{G5.18}. 
\end{proof}

Applying Gronwall's inequality to \eqref{G5.18} and using the inequality \eqref{G5.17},  one has
\begin{align}
\|g\|_{H_{x,v}^2}
\leq\,& e^{CT+C\int_0^t \big(\|\nabla_x u(s)\|_{H^2}^2+\|\nabla_x B(s)\|_{H^2}^2\big)\mathrm{d} s}\Big(\|g_0\|_{H_{x,v}^2}+\delta_3^2 T\Big)\nonumber\\
\leq\,&Ce^{CT}\Big(\|g_0\|_{H_{x,v}^2}+\|(u_0,B_0)\|_{H^2}\Big)\nonumber\\
\leq\,&C (\|g_0\|_{H_{x,v}^2}+\|(u_0,B_0)\|_{H^2})^{\frac{1}{2}}.
\end{align}
Thus, \eqref{N5.3} is justified by selecting $C\tau_1< C \tau_1^{\frac{1}{2}}
<\delta_3$. 


\medskip 
{\it Step 2: Constructing the unique strong solutions to the system \eqref{I8} in $\mathbb{R}^3$.} \ \
Similar to the vacuum case, we construct the strong solutions to the   system \eqref{I8} via 
iteration method. We introduce the approximation sequences $\{(g^{n+1},u^{n+1},B^{n+1})\}$ satisfying:
 \begin{equation}\label{I8a}
 \left\{
\begin{aligned}
& \partial_t g^{n+1}+v \cdot \nabla_x g^{n+1}+\big((v-u^{n}) \times B^{n}\big) \cdot \nabla_v g^{n+1} +(u^{n}\times B^{n})\cdot v M=0, \\
& \partial_t u^{n+1}+u^{n}\cdot \nabla_x u^{n+1}+\nabla_x P^{n+1}-\Delta_x u^{n+1}-B^{n}\cdot\nabla_x B^{n+1}\\
& \qquad  +\frac{1}{2}\sum_{i,j=1}^3 B_{i}^n\partial_j B_i^{n+1}-u^{n+1}\times B^{n}=\int_{\mathbb{R}^3} 
(u^{n+1}\times B^{n}) g^{n+1}\mathrm{d}v-\int_{\mathbb{R}^3} 
(v\times B^{n}) g^{n+1}\mathrm{d}v, \\
& \partial_t B^{n+1}-\Delta_x B^{n+1}+u^{n}\cdot\nabla_x B^{n+1}-B^{n}\cdot\nabla_x u^{n+1}=0, \\
& \nabla_x \cdot u^{n+1}=0, \quad \nabla_x \cdot B^{n+1}=0,
\end{aligned}\right.
\end{equation}
with the initial data: 
\begin{align}\label{I8b}
(g^{n+1}(t, x, v),   u^{n+1}(t, x), B^{n+1}(t, x))|_{t=0}=\,& (g_0(x, v),  u_0(x),  B_0(x)),   
\end{align} 
for all $n\geq 0$ and $(x,v)\in\mathbb{R}^3\times\mathbb{R}^3$. In addition,
\begin{align}\label{I8c}
(g^0(t,x,v),u^0(t,x),B^0(t,x))\equiv(g_0(x,v),u_0(x),B_0(x)),\quad(t,x,v)\in [0,T]\times\mathbb{R}^3\times\mathbb{R}^3.
\end{align} 

{
Based on the above  a priori estimates on the nonlinear system \eqref{I8} obtained in Lemmas 
\ref{L5.3}--\ref{L5.7}, and taking the similar  arguments in Section 2, we obtain that the problem 
\eqref{I8a}--\eqref{I8b} has a unique solution $(g^{n+1},u^{n+1},B^{n+1})$ and 
$(g^{n+1},u^{n+1},B^{n+1})$ converges to $(g,u,B)$ in $C([0,T];H^2_{x,v}(\mathbb{R}^3\times\mathbb{R}^3))\times C([0,T];H^2(\mathbb{R}^3))\times\ C([0,T];H^2(\mathbb{R}^3))$. We omit the details here for brevity. 
In addition, the inequalities \eqref{G51} and \eqref{G5.3} follow from the  above a priori estimates and 
the lower semi-continuity of the norms of $(g^{n+1},u^{n+1},B^{n+1})$. Hence we complete the proof
of Theorem  \ref{T5.2}.}
\end{proof}


\subsection{The periodic domain  $\mathbb{T}^3$ case}
In this subsection,  we study the Vlasov-MHD system \eqref{I8} near
the equilibrium state $(M,0,0)$ in 
  $\mathbb{T}^3$. Assuming that the initial data are sufficient small with 
 an order of $e^{-\mathcal{O}(1)T}$,
we prove the global existence and uniqueness  of strong solutions  on $[0,T]$ for any $T\in (0,+\infty)$  
and provide the exponential decay of $u, B$.  

\begin{proof}[Proof of Theorem \ref{T5.1}]

First, as pointed out in Section 4, due to the boundedness of the periodic domain,  we only need to  assume that $g_0(x,v)$ has a compact support  in the velocity $v$.  We notice that  the existence and uniqueness of strong solutions to \eqref{I8} in  $\mathbb{T}^3$   are essential to those in      $\mathbb{R}^3$   and we omit them for brevity.  
 Below we focus on the uniform a priori estimates to the strong solutions to obtain the exponential decay of $u $ and $B$. 

 By a direct calculation, we have
the following conservation laws:
{\begin{equation}\label{GGG5.21}
\left\{
\begin{aligned}
&\frac{\mathrm{d}}{\mathrm{d} t} \int_{\mathbb{R}^3\times\mathbb{T}^3} g \mathrm{d} v \mathrm{d} x=0,\\
& \frac{\mathrm{d}}{\mathrm{d} t}\left(\int_{\mathbb{T}^3}  u \mathrm{d} x+\int_{\mathbb{R}^3\times\mathbb{T}^3} v g \mathrm{d} v \mathrm{d} x\right)=0, \\
&\frac{\mathrm{d}}{\mathrm{d} t} \int_{\mathbb{T}^3} B \mathrm{d} x=0.
\end{aligned}\right.
\end{equation}}
Under the assumptions of Theorem \ref{T5.1}, we get
 {\begin{equation} \label{GGG5.22}
\left\{
\begin{aligned}
& \int_{\mathbb{R}^3\times\mathbb{T}^3} g \mathrm{d} v \mathrm{d} x=0, \\
& \int_{\mathbb{T}^3}  u \mathrm{d} x+\int_{\mathbb{R}^3\times\mathbb{T}^3} v g \mathrm{d} v \mathrm{d} x=0,  \\  
& \int_{\mathbb{T}^3} B \mathrm{d} x=0 
\end{aligned}\right.
\end{equation}
for all $0 \leq t\leq T$.}

Thanks to the Poincaré's inequality and  \eqref{GGG5.22}, we obtain
\begin{align}\label{G5.24}
\|B\|_{L^2} & \leq C\|\nabla_x B\|_{L^2},\\
\|u\|_{L^2} & \leq C\Big\|\nabla_x u+\int_{\mathbb{R}^3}v\nabla_x g\mathrm{d}v\Big\|_{L^2}+C\Big\|\int_{\mathbb{R}^3}v\nabla_x g\mathrm{d}v\Big\|_{L^2} \nonumber\\
&\leq C\|\nabla_x u\|_{L^2}+C\|g\|_{H_{x,v}^2}\nonumber \\
\label{G5.25}
&\leq C \|\nabla_x u\|_{L^2}+C\tau_2^{\frac{1}{2}} .
\end{align}
Combining \eqref{G5.24} with \eqref{G5.25}, we   arrive at 
\begin{align}\label{G5.26}
\|u\|_{L^2}^2+\|B\|_{L^2}^2\leq C\big(\|\nabla_x u\|_{L^2}^2+ \|\nabla_x B\|_{L^2}^2  \big)  +C\tau_2. 
\end{align}

We define the energy functionals $\mathcal{E}_2(t)$ and the corresponding dissipation rate functional $\mathcal{D}_2(t)$  to the system \eqref{I8} as same as the $\mathbb{R}^3$ case given in \eqref{N5.13}
and \eqref{N5.14}, respectively. 
Similar to the previous arguments on the whole space case,   it holds
\begin{align}\label{G5.27}
\frac{\rm d}{{\rm d}t}\mathcal{E}_{2}(t)+\lambda_{21}\mathcal{D}_{2}(t)\leq 0 
\end{align}
for some $\lambda_{21}>0$.

Therefore, we obtain, from   \eqref{G5.26}, that 
\begin{align}\label{G5.28}
 \frac{\rm d}{{\rm d}t}\mathcal{E}_{2}(t)+\lambda_{22}\mathcal{E}_{2}(t)\leq C\tau_2   
\end{align}
for some $\lambda_{22}>0$.
We define
\begin{align}\label{K5.26}
\tilde{\mathcal{Q}}(t) :=e^{\lambda_{22} t}\mathcal{E}_2(t).
\end{align}
From \eqref{G5.28}, we   achieve that
\begin{align}\label{K5.27}
\frac{\mathrm{d}}{\mathrm{d} t}\tilde{\mathcal{Q}}(t) \leq e^{\lambda_{22} t} \left(\frac{\mathrm{d}}{\mathrm{d} t}\mathcal{E}_2(t)+\lambda_{22}\mathcal{E}_2(t)\right)\leq C e^{\lambda_{22} t}\tau_2.
\end{align}
Thus, 
\begin{align*}
\tilde{\mathcal{Q}}(t) \leq \tilde{\mathcal{Q}(0)}+C\int_0^Te^{\lambda_{22} t}\tau_2     \mathrm{d}t,  
\end{align*}
which implies
\begin{align*}
\mathcal{E}_2(t)&\leq e^{-\lambda_{22} t}\left(\mathcal{E}_2(0)+CTe^{\lambda_{22} T}\tau_2 \right) 
 \leq Ce^{-\lambda_{22} t}(\tau_2+\tau_2^{\frac{1}{2}}) 
 \leq Ce^{-\lambda_{22} t} ,
\end{align*}
where we have utilized the smallness of $\tau_2$. 
Therefore, 
\begin{align*}
& \|( u, B)\|_{H^2} \leq C e^{-\lambda_{23} t}     
\end{align*}
holds for some $\lambda_{23}>0$ and all $0\leq t\leq T$. And we complete the proof of Theorem \ref{T5.1}.
\end{proof}
 
 \smallskip 
 
{\bf Acknowledgements:} 
Li and   Ni  are supported by NSFC (Grant Nos. 12071212, 12331007).  
And Li is also supported by the “333 Project" of Jiangsu Province.

\bibliographystyle{plain}

\end{document}